\documentclass[openright,a4paper,american,11pt]{amsart}

\usepackage[latin1]{inputenc}

\usepackage{amsmath}
\usepackage{amssymb}

\usepackage{babel}

\usepackage{amsfonts}
\input xy
\xyoption{all}
\usepackage[dvips]{graphicx}
\usepackage{boxedminipage}
\usepackage{color}

\newcommand{\ZZ}{{\mathbb Z}}
\newcommand{\NN}{{\mathbb N}}

\newcommand{\CC}{{\mathbb C}}
\newcommand{\RR}{{\mathbb R}}

\newcommand{\m}{{\mathfrak M}}

\newtheorem{thm}{Theorem}[section]
\newtheorem{defi}[thm]{Definition}
\newtheorem{prop}[thm]{Proposition}

          {\theoremstyle{definition}
\newtheorem{rem}[thm]{Remark}}
          {\theoremstyle{definition}
\newtheorem{exa}[thm]{Example}}

\begin{document}
\title{RECURSIVE FORMULAS FOR  WELSCHINGER INVARIANTS of the
  projective plane}
\date{\today}
\author{Aubin Arroyo}
\author{Erwan Brugallé}
\address{Université Pierre et Marie Curie,  Paris 6, 175 rue du Chevaleret, 75 013 Paris, France}
\email{brugalle@math.jussieu.fr}
\author{Lucía López de Medrano}
\address{Unidad Cuernavaca del Instituto de Matemáticas,Universidad Nacional Autonoma de México. Cuernavaca, México}
\email{aubin@matcuer.unam.mx, lucia@matcuer.unam.mx}


\maketitle

\begin{abstract}
 Welschinger invariants of the real projective plane
can be computed via the enumeration of enriched graphs, called marked
floor diagrams. By a purely
combinatorial study of these objects, 
we 
establish
a Caporaso-Harris type formula which allows one to compute Welschinger
invariants 
for configurations of points
with any number of complex 
conjugated points.  

\end{abstract}

\section{Introduction}
\textit{Welschinger invariants} of symplectic 4-manifolds were introduced in
\cite{Wel1b} (see also \cite{Wel1}), and since then they 
are attracting
a
lot of attention. 
One of the 
interests of these invariants is that they
 provide lower bounds in
real enumerative geometry. 
As  an example of application, it is 
proved in \cite{IKS1} that through any 
configuration of $3d-1$ points in $\RR P^2$ passes at least one real
rational algebraic curve of degree $d$ (and even a lot! see also \cite{IKS2}).

Welschinger invariants can be seen as
real analogs of genus 0 \textit{Gromov-Witten invariants}, well known in
complex enumerative geometry. Hence,  
studying  relations between these two sequences of invariants is an
important problem. 
Both
invariants 
can be computed via the same kind of recursive
formulas. In \cite{KonMan1}, Kontsevich 
deduced a recursive formula for
all genus
0 Gromov-Witten invariants of $\CC P^2$  in terms of those of lower
degree
from
the so called \textit{WDVV equation}. Recently, Solomon
\cite{Sol1} 
designated suitable analogues of
WDVV equation in \textit{open Gromov-Witten theory}  
and obtained
similar formulas involving Welschinger 
invariants. Another approach to compute Gromov-Witten invariants
was
proposed by Caporaso and Harris in \cite{CapHar1}. They
obtained
a recursive formula by specializing one point after the other
to lie on a given line in $\CC P^2$. Moreover, this formula computes
not only genus 0 Gromov-Witten invariants but \textit{relative Gromov-Witten
invariants} of any genus and degree. In this paper, we are interested
in a similar formula involving Welschinger invariants.

\vspace{1ex}
Thanks to Mikhalkin's Correspondence Theorem (see \cite{Mik1}),
\textit{tropical geometry} turned out to be a powerful tool to solve a lot of 
enumerative problems (see for example \cite{Mik1}, \cite{Sh8},\cite{IKS1}, \cite{Br6b},
\cite{Br7}). In particular, 
it provides a quite combinatorial and very practical way to compute
Gromov-Witten and Welschinger 
invariants.  

In \cite{GM2}, Gathmann and Markwig 
gave a tropical version and a tropical proof of Caporaso-Harris formula.
Then, 
  Itenberg, Kharlamov and Shustin 
adapted 
in \cite{IKS3} 
this tropical approach to a real setting and produced a
Caporaso-Harris like recursive formula for Welschinger invariants in
the totally real case. 
Their formula
 involves not only Welschinger invariants, but also other numbers
  which reveal to be tropical relative Welschinger invariants.
This
  was quite surprising since 
in a non-tropical setting
no relative Welschinger invariants
are 
known.

An even more combinatorial way to compute Gromov-Witten and
Welschinger invariants has been proposed
in \cite{Br6b} and \cite{Br7}, where tropical curves are replaced by
\textit{floor diagrams}. 
Gathmann-Markwig tropical formula in the complex case and 
 Itenberg-Kharlamov-Shustin's one in the totally real case are easily
 translated into the language of floor diagrams. 
Moreover, the quite simple
combinatoric of these objects allows one to generalize easily the
formula in \cite{IKS3} to Welschinger invariants for 
collections of points mixing real points and pairs of complex
conjugated ones.

\vspace{1ex}
In this paper, we explain how to pass from marked floor diagrams to 
a Caporaso-Harris style  formula
which allows one 
to compute all Welschinger invariants of the projective plane. 
This
work is based on \cite[Theorem 4.16]{Br6b}, which in its turn is
based on the two Correspondence Theorems by Mikhalkin (\cite{Mik1}) and 
Shustin (\cite{Sh8}).

Our
formula does not involve directly Welschinger invariants, but some
other numbers, the sum of which gives Welschinger invariants. In
regard to tropical relative Welschinger invariants for configurations of real
points, it is natural to wonder about the invariance of the numbers
involved in our formula. However, one sees easily that they do not
correspond to any invariant, even tropical. Worse,  we give an
example in section \ref{trop rel inv} which shows that there is no hope
to define as in \cite{IKS3} tropical relative Welschinger invariants
for real configurations of points containing some pairs of complex
conjugated points.

\vspace{2ex}
Another approach for computing Welschinger invariants, based on the
symplectic field theory, has been proposed by Welschinger in
\cite{Wel4}.
There he derived closed graph-combinatorial
formulas in certain cases, in particular for the real projective 
plane and real configurations consisting of 
pairs of complex conjugated points together with up to two real points.

\vspace{2ex}

The structure of the paper is the following: In section \ref{conv} we fix notation and convention used thereafter. In section \ref{enum} we define Gromov-Witten and Welschinger invariants, and in section \ref{floor} we explain how to compute these enumerative invariants using floor
diagrams. Section \ref{formula} is devoted to the statement of our
formula, and the proof is given in section \ref{proof}.
Some final remarks and computations are given in the last section.


\vspace{2ex}

\textbf{Acknowledgement : }Part of this work was done during the stay
of the last two authors at the Centre Interfacultaire Bernoulli in
Lausanne. We thank the CIB and the National Science Foundation for
excellent working conditions. First and third authors were partially supported by CONACyT 58354 and
37035, and UNAM-PAPIIT IN102307 and IN105806.

We are also grateful to the anonymous referees for their valuable remarks and
suggestions which helped us to improve the initial text.
\section{Convention}\label{conv}

We fix here notation and  convention we use in
this paper.  

\subsection{Notation}
We denote by $\NN$ the set of nonnegative integer numbers and  by
$\NN^*$ the set of positive integer numbers, i.e. $\NN=\{n\in\ZZ 
| \ n\ge 0\}$ and $\NN^*=\{n\in\ZZ 
| \ n> 0\}$.
The set of sequences of elements in $\NN$  having only finitely 
many non
zero terms is denoted by $\NN^\infty$. We denote by $e_i$ the vector  in
$\NN^\infty$ whose all
coordinates are 0 but the 
$i^{th}$ which is equal to 1.
If $\alpha$ is a vector in $\NN^\infty$, we denote by $(\alpha)_i$
its $i^{th}$ coordinate. Given two vectors $\alpha$ and $\beta$ in
$\NN^\infty$, we write
$\alpha\ge \beta$ if $(\alpha)_i\ge (\beta)_i$ for all $i$.
Finally, we put 
$$|\alpha|=\sum_{i=1}^\infty (\alpha)_i, \ \ \ \ \ \
I\alpha = \sum_{i=1}^\infty i(\alpha)_i, \ \ \ \ \ \
I^\alpha = \prod_{i=1}^\infty i^{(\alpha)_i} $$

If $a$ and $b$ are two integer numbers, $\left(\begin{array}{c}
a\\ b\end{array} \right) $ denotes the binomial coefficient (i.e. is
equal to $\frac{a!}{b!(a-b)!}$ if $0\le b\le a$ and to $0$
otherwise). If $a$ and $b_1,b_2,\ldots,b_k$ are 
integer numbers then $\left(\begin{array}{c}
a
\\ b_1,\ldots, b_k
\end{array} \right) $ denotes the multinomial coefficient, i.e.
$$\left(\begin{array}{c}
a
\\ b_1,\ldots, b_k
\end{array} \right) = \prod_{i=1}^{k}\left(\begin{array}{c}
a -\sum_{j=1}^{i-1}b_j
\\ b_i
\end{array} \right)$$


If $\alpha, \alpha_1,\ldots, \alpha_l$ are vectors in $\NN^\infty$, then we
put

$$\left(\begin{array}{c}
\alpha
\\ \alpha_1,\ldots, \alpha_l
\end{array} \right) = \prod_{i=1}^\infty
 \left(\begin{array}{c}
(\alpha)_i
\\ (\alpha_1)_i,\ldots, (\alpha_l)_i
\end{array} \right) $$

\subsection{Graphs}

An 
\textit{oriented graph} 
$\Gamma$ is a 
pair $ (
\overline{\text{Vert}}(\Gamma) ,\text{Edge}(\Gamma) )$  
of two finite sets
$\overline{\text{Vert}}(\Gamma)$ and $\text{Edge}(\Gamma)$, where
$\text{Edge}(\Gamma)$ is a 
list of elements 
of
$\overline{\text{Vert}}(\Gamma)\times
\overline{\text{Vert}}(\Gamma)$.
Note that $\text{Edge}(\Gamma)$ might contain repeated elements.
 An element of
$\overline{\text{Vert}}(\Gamma)$ is called a \textit{vertex} of
$\Gamma$. An element $(v_1,v_2)$ of $\text{Edge}(\Gamma)$ is called an
\textit{edge} of $\Gamma$ and is said to be oriented from  $v_1$ to $v_2$.
A vertex $v$ and an edge $e$ are said to be \textit{adjacent} if $e=(v,v')$ or
$(v',v)$.
A vertex of $\Gamma$ such that all its adjacent edges are outgoing is
called a \textit{source}.
We denote by
$\text{Vert}^{\infty}(\Gamma)$ the set of sources of
$\Gamma$, and we put
$\text{Vert}(\Gamma)=\overline{\text{Vert}}(\Gamma)\setminus 
\text{Vert}^{\infty}(\Gamma)$. We also denote by 
$\text{Edge}^{\infty}(\Gamma)$ the set of edges adjacent to a
source.
A 
graph
$\Gamma$ is  \textit{connected} if for any two elements 
 $v_1$ and $v_2$ of $\overline{\text{Vert}}(\Gamma)$, there exists a sequence
 $s_1=v_1$, $s_2$, $\ldots$, $s_{k-1}$,
$s_k=v_2$ of elements of $\overline{\text{Vert}}(\Gamma)$ such that
$(s_i,s_{i+1})\in\text{Edge}(\Gamma)$ or 
$(s_{i+1},s_{i})\in\text{Edge}(\Gamma)$ for all $i$. 
A \textit{tree} is a connected graph with Euler characteristic equal
to 1 (i.e. $\sharp \overline{\text{Vert}}(\Gamma) - \sharp
\text{Edge}(\Gamma)  =1$). In this case, $\text{Edge}(\Gamma)$ is
actually a subset of $\overline{\text{Vert}}(\Gamma)\times
\overline{\text{Vert}}(\Gamma)$.

An oriented tree $\Gamma$ is naturally enhanced with a partial ordering : an
element 
(i.e. an edge or a vertex of $\Gamma$)
$a$ of $\Gamma$ is greater than another element  $b$ if there
exists a sequence $c_1=a$, $c_2$, $\ldots$, $c_{k-1}$,
$c_k=b$ of elements of $\Gamma$ such that $c_i$ is adjacent to
$c_{i+1}$ for all $i$, and if $c_i$ (resp. $c_{i+1}$) is an edge of $\Gamma$ then
$c_{i}=(v,c_{i+1})$ (resp. $c_{i+1}=(c_i,v)$).   

A \textit{weighted graph} is a graph $\Gamma$ enhanced with 
  a function
$\omega:\text{Edge}(\Gamma)\rightarrow \mathbb N^*$. The integer $\omega(e)$ is called the \textit{weight}
  of the edge $e$.
 The weight allows
one to define the \textit{divergence} at the vertices. Namely, for a
vertex $v\in\overline{\text{Vert}}(\Gamma)$ we define the divergence
$\text{div}(v)$ to be the sum of the weights of all incoming edges
minus the sum of the weights of all outgoing edges.

 


\section{Complex and real enumerative problems}\label{enum}

\subsection{Relative Gromov Witten invariants}

Fix $d\ge 1$ an integer number,  $\omega=\{p_1,\ldots, p_{3d-1}\}$ a
generic collection of $3d-1$ points in $\CC P^2$ and denote by $\mathcal
C(\omega)$ the set of all 
 irreducible complex
rational curves of degree $d$ in $\CC P^2$ containing
$\omega$. It is well known that the cardinal of $\mathcal
C(\omega)$ does not depend on $\omega$, and it is called  
the genus 0 \textit{Gromov-Witten
invariant} of degree $d$  of  $\CC P^2$.

More generally, one can define \textit{relative Gromov-Witten
  invariants}. Fix  a line $L$ in $\CC
P^2$,
 a degree $d\ge 1$ and two
vectors $\alpha$ and $\beta$ in $\NN^\infty$ such that $I\alpha
+I\beta=d$. Choose $\omega=\{p_1,\ldots, p_{2d-1+|\beta|}\}$ a
collection of  $2d-1+|\beta|$ points in $\CC P^2$, and
$\omega_L=\{p_1^1,\ldots,$
$ p_{(\alpha)_1}^1,p_1^2,\ldots,
p_{(\alpha)_2}^2, \ldots,$
$p_1^k,\ldots,
p_{(\alpha)_k}^k, \ldots\}$  a collection of $|\alpha|$ points on
$L$. Denote by $\mathcal
C(\omega,\alpha,\beta)$ the set of all 
 irreducible complex
 rational curves of degree $d$ in $\CC P^2$ containing
$\omega$, 
having no singular point on $L$,
intersecting the line $L$ 
at the point $p_i^j\in \omega_L$
with multiplicity $j$ for all $i$ and $j$, and intersecting $L$ at $(\beta)_j$ 
additional points with multiplicity $j$ for all $j$. 
Once again, the cardinal of $\mathcal
C(\omega,\alpha,\beta)$ does not depend on $\omega$ and $\omega_L$ as
long as these configurations are generic, and is denoted by
$N^{\alpha,\beta}(d)$. The number $N^{0,de_1}(d)$ is the 
genus 0 Gromov-Witten invariant of degree $d$.


Note that one can define relative Gromov-Witten invariants of any
genus, and that the  Caporaso-Harris formula computes all these
numbers. However, we are only interested in rational curves in this paper, and
the formula we write in section \ref{ch} only deals with rational curves.

\subsection{Welschinger invariants}\label{real}
Let us now 
consider
real algebraic curves.
Fix 
 an integer number $d\ge 1$,  and  a
generic collection $\omega=\{p_1,\ldots, p_{3d-1}\}$
of $3d-1$ points in $\CC P^2$.
Suppose 
that $\omega$ is real, i.e. $conj(\omega)=\omega$ where $conj$ is the
standard complex conjugation on $\CC P^2$.
Consider the set $\RR \mathcal
C(\omega)$  of all 
 irreducible real 
 rational curves of degree $d$ in $\CC P^2$ containing
$\omega$. The cardinal of $\RR \mathcal
C(\omega)$ is no longer invariant with respect to $\omega$, but
Welschinger proved in \cite{Wel1} that if one 
counts the curves in 
 $\RR \mathcal
C(\omega)$ with an appropriate sign, then one  obtains an
invariant. 

If $C$ is a real 
 algebraic 
nodal curve in $\CC P^2$, let us denote by
$w(C)$ the number of isolated double points of $C$ in $\RR P^2$
(i.e. points where $C$ has local equation $X^2+Y^2=0$). 

\begin{thm}[Welschinger]
The number
$$\sum_{C\in \RR \mathcal
C(\omega)}  (-1)^{w(C)}$$
only depends on the degree $d$ and the number of pairs of complex
conjugated points in $\omega$.
\end{thm}
These numbers are called \textit{Welschinger invariants}, and 
we denoted them by $W_2(d,r)$ where $r$ is the number of pairs of complex
conjugated points in $\omega$.




Of course, one can count with Welschinger
signs real algebraic curves of any genus or with tangency
conditions. However the number obtained is no longer an invariant (see
\cite{Wel1} or \cite{IKS3}). See section \ref{trop rel inv} for this
discussion in the tropical setting.

\section{Floor diagrams}\label{floor}

Here we define floor diagrams and state their relation with
Gromov-Witten and Welschinger invariants of $\CC P^2$. The definitions
of this paper differ slightly  
from those of \cite{Br6b} and
\cite{Br7}. The first reason is that here we are mainly interested in
enumeration 
of plane curves,  for which  
the 
floor diagrams we have
to consider are 
much simpler than those required in arbitrary
dimension. 
On the other hand, 
the
floor diagrams which are needed to
compute absolute invariants are not sufficient to deal with relative
invariants, we have to allow edges in $\text{Edge}^\infty(\mathcal
D)$ with any positive weight.

\begin{defi}\label{def fd}
 A connected weighted oriented tree $\mathcal{D}$ is
  called a \text{floor diagram} of genus $0$ and degree $d$ if the
  following conditions hold 
\begin{itemize}
\item for any $v\in\text{Vert}(\mathcal{D})$, one has $\text{div}(v)=1$,
\item 
  each source has a unique adjacent edge, 
\item one has $\sum_{v\in\text{Vert}^{\infty}(\mathcal{D}) }div(v)=-d$.
\end{itemize}
\end{defi}
Note that Definition \ref{def fd} implies that the set $\text{Vert}(\mathcal
D)$ has exactly 
$d$ elements.

\begin{figure}[h]
\centering
\begin{tabular}{cccccccc}
\includegraphics[height=3cm, angle=0]{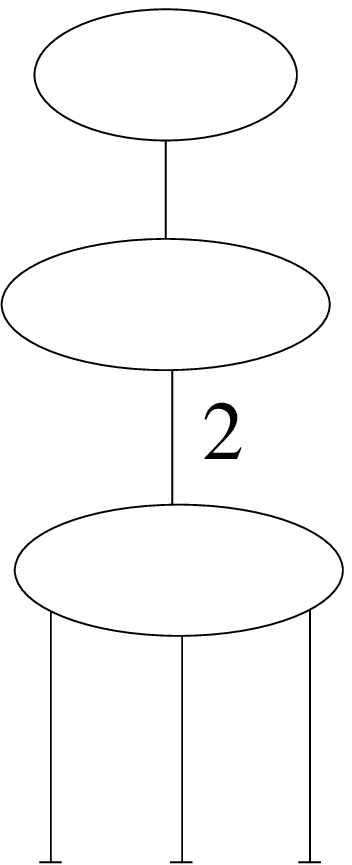}&
\includegraphics[height=3cm, angle=0]{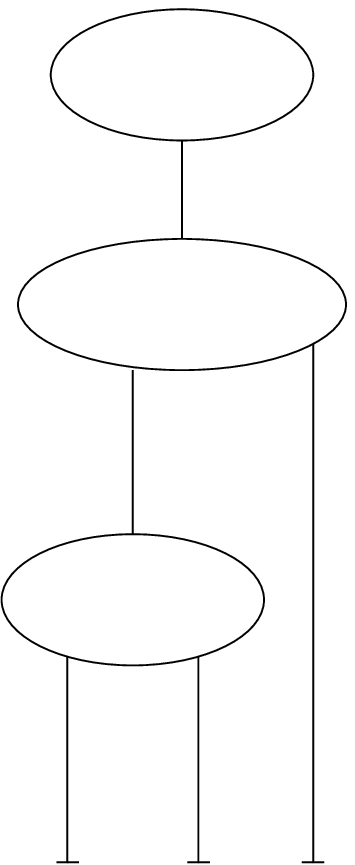}&
\includegraphics[height=3cm, angle=0]{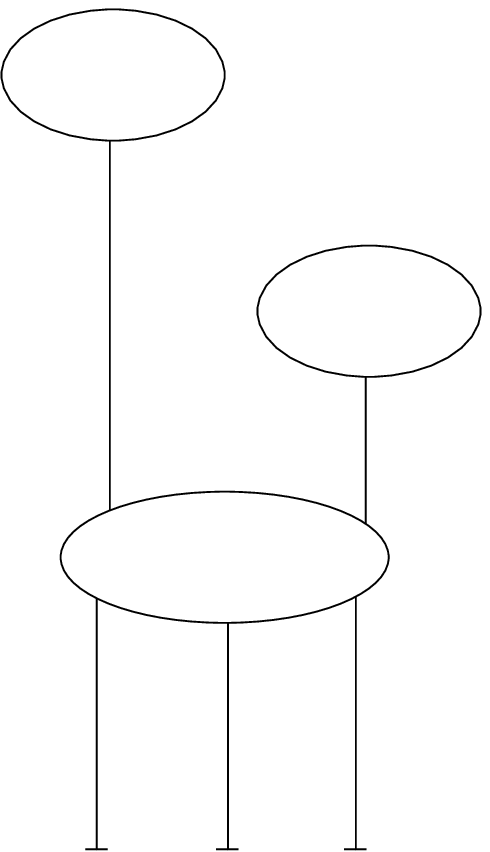}&
\includegraphics[height=3cm, angle=0]{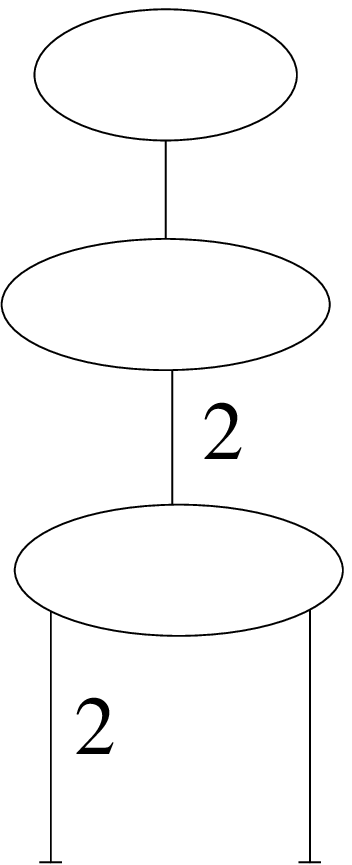}&
\includegraphics[height=3cm, angle=0]{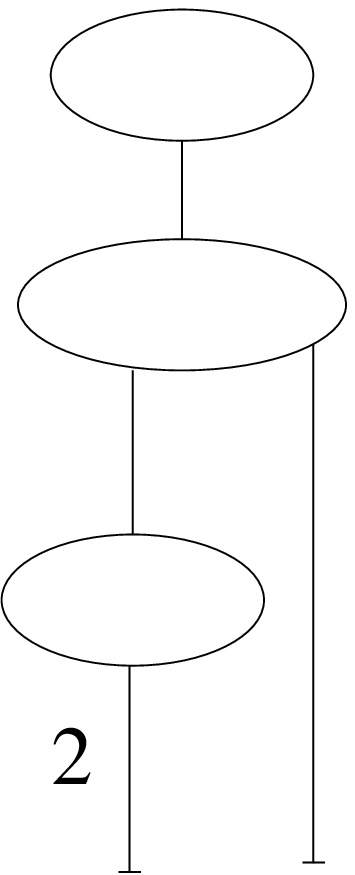}&
\includegraphics[height=3cm, angle=0]{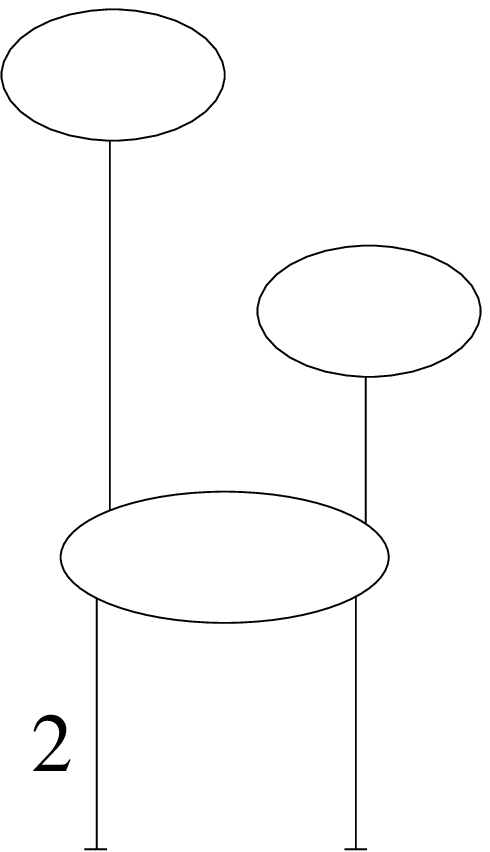}&
\includegraphics[height=3cm, angle=0]{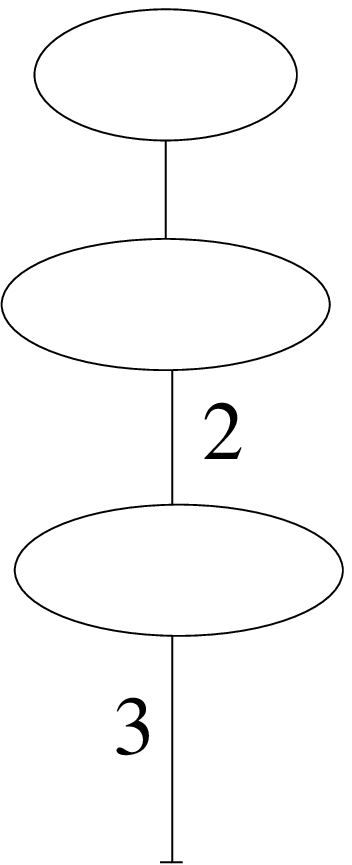}&
\includegraphics[height=3cm, angle=0]{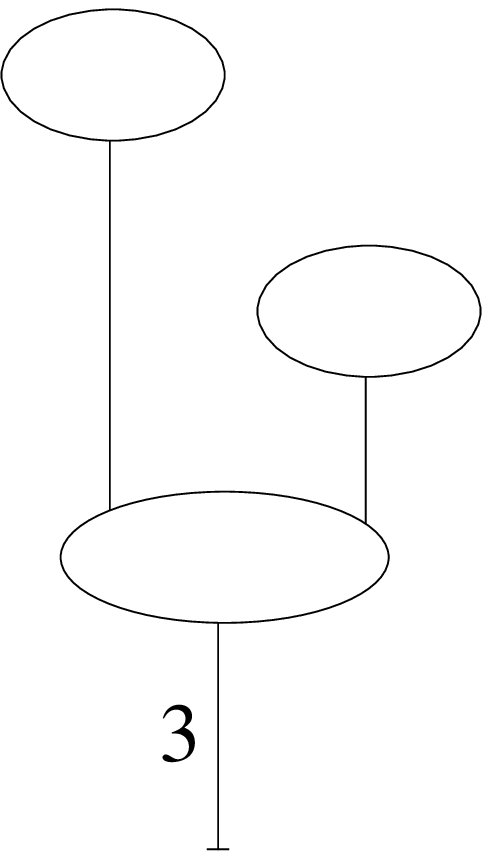}
\end{tabular}
\caption{Floor diagrams of degree 3 and genus 0}
\label{FD3}
\end{figure}

Here are the convention we use to depict 
floor diagrams : 
elements of 
$\text{Vert}^\infty(\mathcal D)$
are represented by small horizontal segments,
elements of 
$\text{Vert}(\mathcal D)$
are represented by ellipses. Elements of 
$\text{Edge}(\mathcal D)$ 
are represented by vertical lines, and 
the orientation is
implicitly from down to up. The weight of an edge 
is
indicated 
only 
when
it is 
at least 2.
All floor diagrams of degree 3
and genus 0 are depicted in Figure \ref{FD3}.

Let $\alpha,\beta,\gamma,\delta$ be four vectors in $\NN^\infty$, and define
$n=2d-1+|\alpha+\beta+2\gamma+2\delta|$ and
$\mathcal{P}=\{1,\dots,n\}$. Let $\m :\mathcal{P}\rightarrow
\mathcal{D}$ be a map. 

\begin{defi}\label{mfd}
The map $\m$ is called a marking of $\mathcal{D}$  of type 
    $(\alpha,\beta,\gamma,\delta)$  if the following
  conditions hold
\begin{itemize}
\item the floor diagram $\mathcal{D}$ is  of genus $0$ and of degree
  $I\alpha+I\beta+2I\gamma+2I\delta$, 
\item the map $\m$ is injective,
\item if $\m(i)>\m(j)$, then $i>j$,
\item if $\sum_{j=1}^{k-1}(\alpha)_j+1  \le i\le \sum_{j=1}^{k}(\alpha)_j$
  or  $|\alpha|+2\sum_{j=1}^{k-1}(\gamma)_j+1  \le i\le
  |\alpha|+2\sum_{j=1}^{k}(\gamma)_j$, then
  $\m(i)$ is a source with divergence $-k$, 

\item for any $k \ge 1$, there are exactly $(\beta)_k+2(\delta)_k$ elements of
  $\text{Edge}^{\infty}(\mathcal{D})$ with weight $k$   in the
  image of $\m$, 

\item for any source $v$ of $\mathcal D$ adjacent to the edge $e$, 
  exactly one of the two elements $v$ or $e$ is in the image of $\m$.

\end{itemize}
\end{defi}
A floor diagram enhanced with a marking of type
$(\alpha,\beta,\gamma,\delta)$ is called a \textit{marked floor
  diagram of type $(\alpha,\beta,\gamma,\delta)$} and is said to be
marked by $\m$.
Two marked floor diagrams are called \textit{equivalent} if they can
be 
identified by a homeomorphism of oriented graphs. 

A simple Euler characteristic
computation shows that
if $\m$ is a marking of $\mathcal D$, then any vertex in $\text{Vert}(\mathcal D)$
and any edge in
$\text{Edge}(\mathcal D)\setminus \text{Edge}^\infty(\mathcal D)$
is  in the image of $\m$.
To any marked floor diagram, we  assign a sequence of nonnegative
numbers called \textit{multiplicities}~: a \textit{complex}
multiplicity, and some \textit{$r$-real} multiplicities. 

\begin{defi}
The complex multiplicity of a marked floor diagram $\mathcal D$ of type
$(\alpha,\beta,\gamma,\delta)$, denoted by $\mu^\CC(\mathcal D)$, is defined as
$$\mu^\CC(\mathcal D)=I^{-2\alpha-\beta-4\gamma-2\delta}
\prod_{e\in \text{Edge}(\mathcal D)}w(e)^2 $$  
\end{defi}
Note that the complex multiplicity of a marked floor diagram 
depends only on the underlying floor diagram.
Our formula mixes real and complex invariants (see
section \ref{fW}), hence we will need the following theorem in section
\ref{proof}.  
This theorem 
is proved in \cite{Br6b} (see also \cite{Br7}) in
the case $\alpha=0$ and $\beta=de_1$. However, it is straightforward
that the
same proof works for any couple $(\alpha,\beta)$ (see also \cite{GM2}).

\begin{thm}[Brugallé, Mikhalkin]\label{NFD}
For any $\alpha$ and $\beta$ in $\NN^\infty$, and $d=I\alpha+I\beta$, one has
$$N^{\alpha,\beta}(d)=\sum \mu^\CC(\mathcal D)$$
where the sum is taken over all marked floor diagrams of degree $d$,
genus $0$ and of type $(\alpha,\beta,0,0)$.
\end{thm}

\begin{rem}
Taking the sum over  all  marked floor diagrams of degree $d$,
genus $0$ and of type $(\alpha,\beta,\gamma,\delta)$, one obtains the
number $N^{\alpha+2\gamma,\beta+2\delta}(d)$.
\end{rem}


\vspace{2ex}

Let us fix $r\ge 0$ such that
$2d-1+|\beta+2\delta|-2r\ge 0$, and $\mathcal D$ a
floor diagram of type $(\alpha,\beta,\gamma,\delta)$ marked by $\m$.

The set $\{i,i+1\}$ is a called \textit{$r$-pair} if 
$i=|\alpha| +  2k -1 $ with 
 $1\leq k \leq |\gamma| $
or
 $i=n-2k+1$ with  $1\leq k \leq r$.
 Denote by
$\Im (\mathcal{D},\m,r)$ the union of all the sets $\{\m(i),
\m(i+1)\}$ where $\m(i)$ is not adjacent to  $\m(i+1)$ and $\{i,i+1\}$ is an 
$r$-pair.
Let $\rho_{\mathcal D,\m,r}:\mathcal{D}\rightarrow \mathcal{D}$ be
the bijection 
defined by  $\rho_{\mathcal D,\m,r} (a)=a$ if $a\in
\mathcal D\setminus \Im(\mathcal D,\m,r)$, and by $\rho_{\mathcal
  D,\m,r}(\m(i))=\m(j)$ if 
$\{i,j\}$ is an $r$-pair  and $\{\m(i),\m(j)\}\subset \Im(\mathcal
D,\m,r)$. Note that  $\rho_{\mathcal D,\m,r}$ is an involution. 


\begin{defi}\label{defi real}
A marked floor diagram $\mathcal D$ of type  $(\alpha,\beta,\gamma,\delta)$
 is called $r$-real if
$(\mathcal{D},\m)$ and
$(\mathcal{D},\rho_{\mathcal D,\m,r}\circ\m)$ are
equivalent, and if exactly $2(\delta)_k$ edges of weight $k$   are in 
$\text{Edge}^{\infty}(\mathcal{D})\cap\Im(\mathcal D,\m,r)$
for any $k\ge 1$. 

The $r$-real multiplicity of a  marked floor
diagram of type  $(\alpha,\beta,\gamma,\delta)$, denoted by
$\mu^\RR_r(\mathcal D,\m)$, is defined as 
$$\mu^\RR_r(\mathcal D,\m)= (-1)^{\frac{\sharp (\text{Vert}(\mathcal D) \cap
  \Im(\mathcal D,\m,r) )}{2}} I^{ -\delta} 
\prod_{
e\in \text{Edge}(\mathcal D) \setminus  \m(\{1,\ldots, n-2r\})}w(e) $$
if $(\mathcal D,\m)$ is an $r$-real marked floor diagram with all
edges of even weight in $\Im(\mathcal D,\m,r)$, and as
$$\mu^\RR_r(\mathcal D,\m)=0 $$
otherwise. 

\end{defi}

Note that as soon as $r\ge 1$, the $r$-real multiplicity of an
$r$-real marked floor 
diagram depends not only on the underlying floor diagram but also on
the marking. 
The following theorem is proved in \cite{Br6b}.
\begin{thm}[Brugallé, Mikhalkin]\label{WFD}
For any $d\ge 1$ and any $r\ge 0$ such that $3d-1-2r\ge 0$, one has
$$W_2(d,r)=\sum \mu^\RR_r(\mathcal D,\m)$$
where the sum is taken over all $r$-real marked floor diagrams of degree $d$,
genus $0$ and of type $(0,(d-2i)e_1,0,ie_1)$ with $0\le i\le \frac{d}{2}$.
\end{thm}

\begin{exa}
All marked floor diagrams of degree 3, genus 0 and type
$(0,(3-2i)e_1,0,ie_1)$  are depicted in 
Table \ref{cubic} with their multiplicities. The first floor diagram
has an edge of weight 2, but we didn't mention it in the picture to
avoid confusion.
According to Theorems
\ref{NFD} and \ref{WFD} one finds $N^{0,3e_1}(3)=12$ and $W_2(3,r)=8-2r$.
\end{exa}

\begin{table}[h]
\centering
\begin{tabular}{c|c|c|c|c|c|c|c|c|c}
&
\includegraphics[height=2.5cm, angle=0]{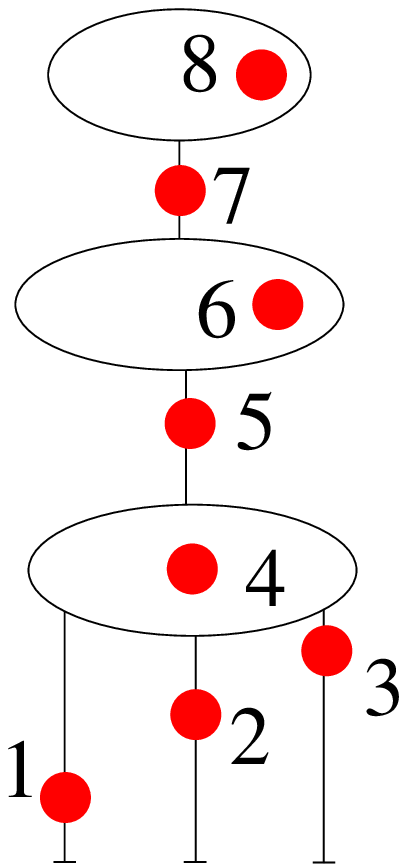}&
\includegraphics[height=2.5cm, angle=0]{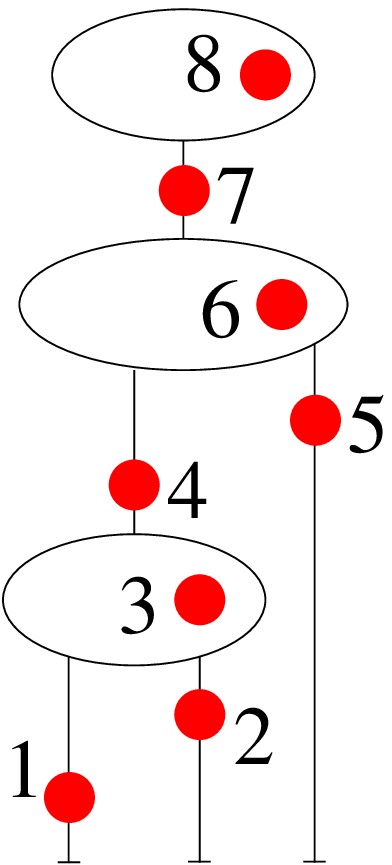}&
\includegraphics[height=2.5cm, angle=0]{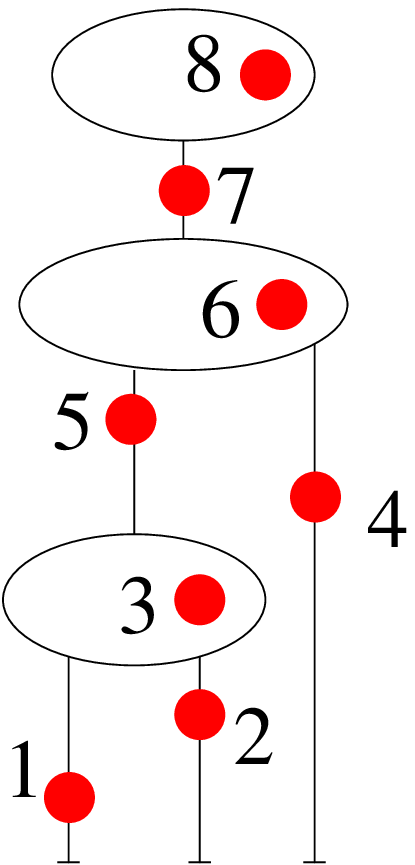}&
\includegraphics[height=2.5cm, angle=0]{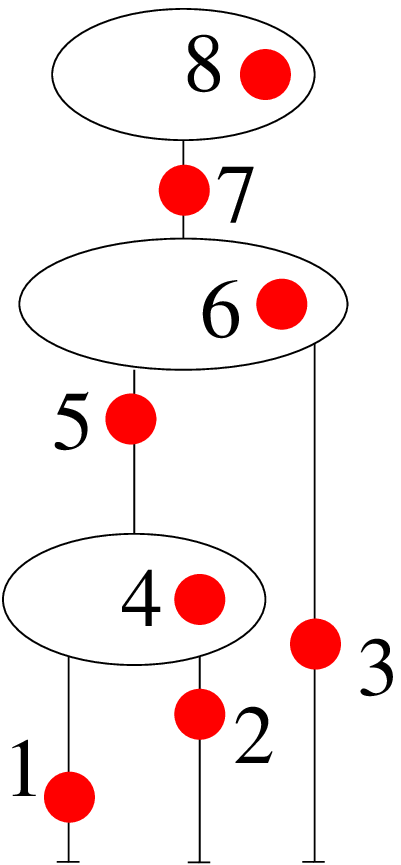}&
\includegraphics[height=2.5cm, angle=0]{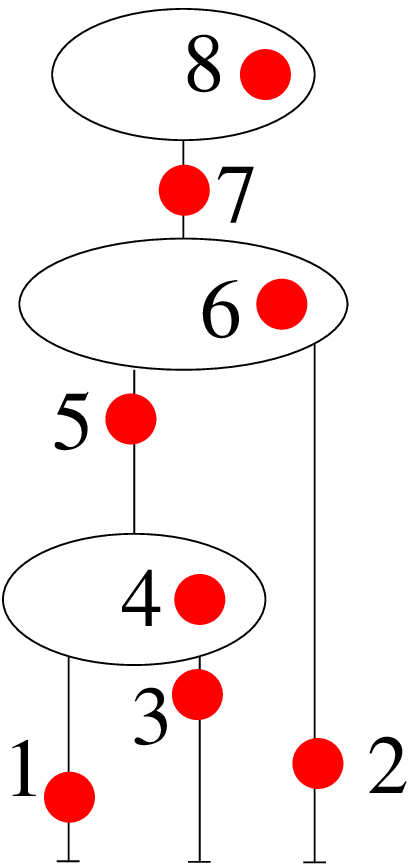}&
\includegraphics[height=2.5cm, angle=0]{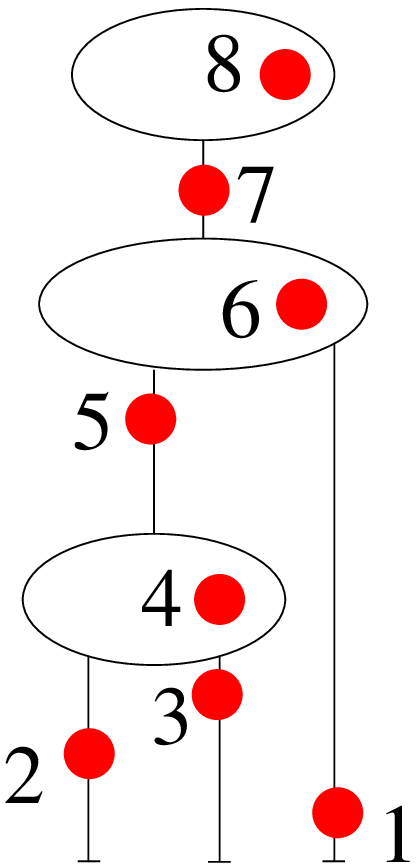}&
\includegraphics[height=2.5cm, angle=0]{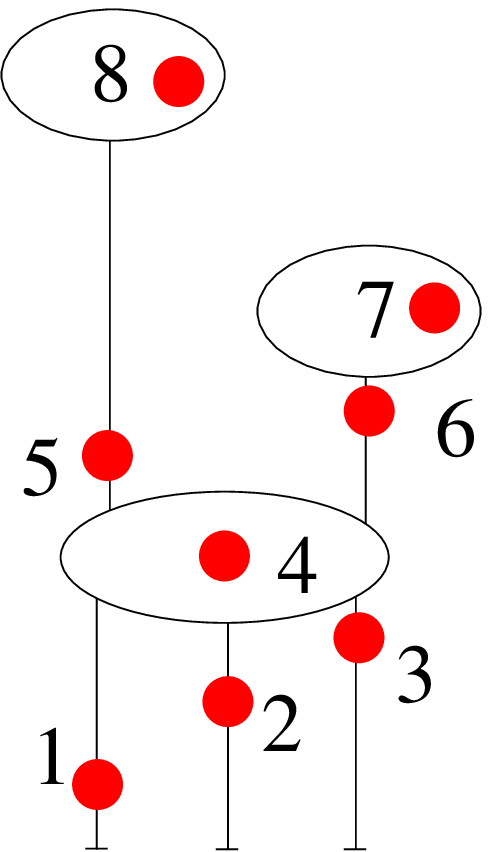}&
\includegraphics[height=2.5cm, angle=0]{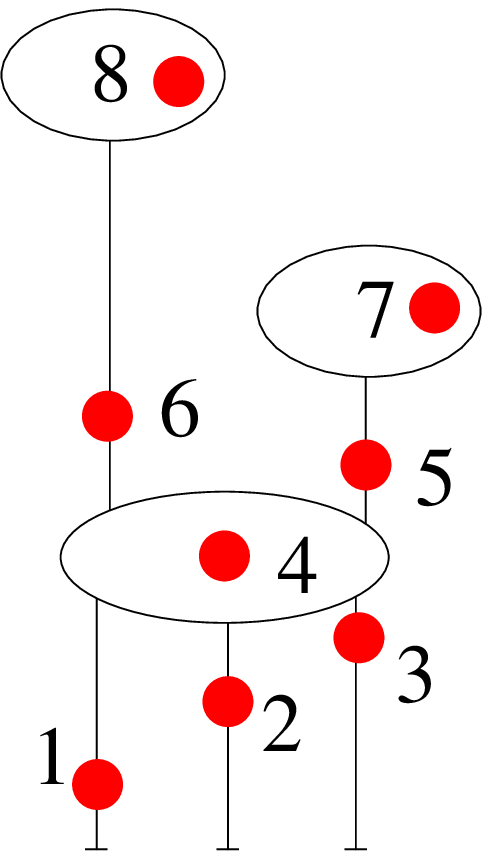}&
\includegraphics[height=2.5cm, angle=0]{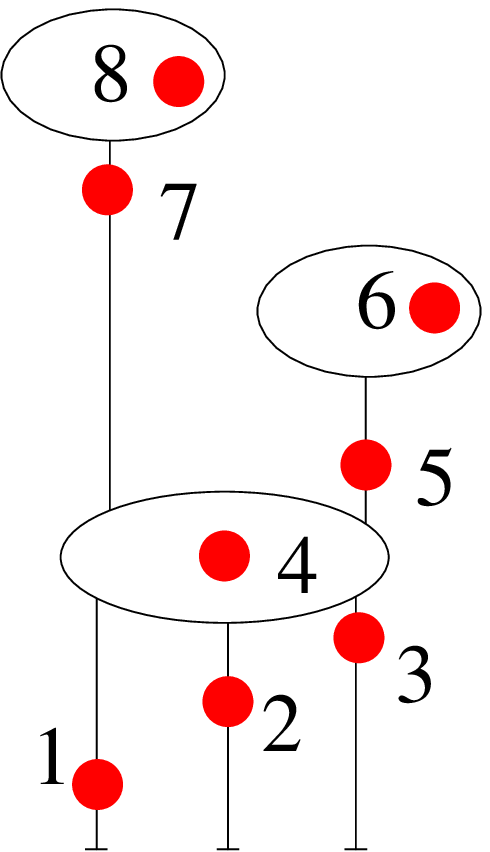}
\\ \hline $ \mu^\CC$ & 4 & 1 &1 &1 &1 &1 &1 &1 &1 
\\ \hline $ \mu^\RR_0$ & 0 & 1 &1 &1 &1 &1 &1 &1 &1 
\\ \hline $ \mu^\RR_1$ & 0 & 1 &1 &1 &1 &1 &0 &0 &1 
\\ \hline $ \mu^\RR_2$ & 0 & 1 &1 &1 &1 &1 &-1 &-1 &1 
\\ \hline $ \mu^\RR_3$ & 0 & 1 &0 &0 &1 &1 &-1 &-1 &1
\\ \hline $ \mu^\RR_4$ & 0 & 1 &0 &0 &0 &0 &-1 &-1 &1 

\end{tabular}
\\
\begin{tabular}{c}
\end{tabular}
\caption{Marked floor diagrams of degree 3 and genus 0}
\label{cubic}
\end{table}

\section{Recursive formulas}\label{formula}

\subsection{Caporaso-Harris formula}\label{ch}
As a warm up for Theorem \ref{form}, we first remind the
Caporaso-Harris formula. Note that the equivalence relation $\sim_s$
will be used in section \ref{fW}.

Given  two integer numbers $l\ge 0$ and $d\ge 0$ and two vectors
$\alpha$ and $\beta$ 
in $\NN^\infty$, we denote by $\mathcal S(d,l,\alpha,\beta)$ the set
composed by the vectors
$(d_1,\ldots,d_l,$
$k_1,\ldots,k_l,$
$\alpha_1,\ldots,\alpha_l,$
$\beta_1,
\ldots,\beta_l)\in   
(\NN^*)^{2l}\times (\NN^\infty)^{2l}$
satisfying 

\begin{itemize}
\item $\forall i, (d_i,k_i,\alpha_i,\beta_i)\le
  (d_{i+1},k_{i+1},\alpha_{i+1},\beta_{i+1})$ for the lexicographic order, 
\item $\sum d_i = d-1$,
\item $\sum \alpha_i \le \alpha$,
\item $\forall i, \beta_i\ge e_{k_i}$, 
\item $\sum (\beta_i-e_{k_i})
 =\beta $,

\item $\forall i, I\alpha_i +I\beta_i=d_i$.

\end{itemize}

To any element 
$s=(
d_1,\ldots,d_l,k_1,\ldots,k_l,\alpha_1,\ldots,\alpha_l,
\beta_1,\ldots,\beta_l)$     
of $\mathcal S(d,l,\alpha,\beta)$, we associate the equivalence
relation $\sim_s$ 
on the set $\{1,\ldots ,l\}$ defined by

$$i\sim_s j \Leftrightarrow (d_i,k_i,\alpha_i,\beta_i)=
(d_j,k_j,\alpha_j,\beta_j)$$ 
For each of the equivalent classes of $\sim_s$, evaluate the factorial of its
cardinal, and denote by $\sigma(s)$ the product of these factorials.

\begin{thm}[Caporaso, Harris]\label{CH}

The numbers $N^{\alpha,\beta}(d)$ are given by the initial value
$N^{e_1,0}(1)=1$ and the relation

\vspace{2ex}

$\begin{array}{lcl}
N^{\alpha ,\beta}(d) &=&\sum_{k|\beta\ge e_k}kN^{\alpha+e_k,
  \beta-e_k}(d) \ \ +
\\
\\ && \sum_{\begin{array}{l}_{l\ge 0}\\ _{s\in \mathcal
  S(d,l,\alpha,\beta)}\end{array}} \left[\frac{1}{\sigma(s)}\left(\begin{array}{c} 
2d-2 +|\beta|
\\ 2d_1 -1+|\beta_1|,\ldots, 2d_l -1+|\beta_l|
\end{array} \right) \right.
\\
\\ &&\left. \left(\begin{array}{c}
\alpha
\\ \alpha_1,\ldots, \alpha_l
\end{array} \right) 
\prod_{i=1}^l (\beta_i)_{k_i}k_i N^{\alpha_i,\beta_i}(d_i) \right]

\end{array}
$

\end{thm}

\subsection{Formula for plane Welschinger invariants}\label{fW}

We state now the main result of this paper, a recursive formula in the
spirit of Theorem \ref{CH} which allows one to compute the numbers
$W_2(d,r)$. Let  $\alpha$, $\beta$,
$\gamma$ and $\delta$ be four vectors in $\NN^\infty$, define
$d= I\alpha  +I\beta +2I\gamma +2I\delta$, and choose an
integer number $r\ge 0$  such
that $2d-1 +|\beta+2\delta|-2r\ge 0$. 
Our recursive formula does not involve the numbers $W_2(d,r)$, but the
numbers $C^{\alpha,\beta,\gamma,\delta}(d,r)$ which are defined as follows

$$C^{\alpha,\beta,\gamma,\delta}(d,r) = \sum \mu^\RR_r(\mathcal
D,\m)$$
where the sum is taken over all marked floor diagrams of type
$(\alpha,\beta,\gamma,\delta)$.
Note that only $r$-real  marked floor diagrams contribute to this sum.
For any $d\ge 1$ and $0\le r\le  \frac{3d-1}{2}  $, 
Theorem \ref{WFD} states that

$$W_2(d,r)=\sum_{i=0}^{\frac{d}{2}}  C^{0,(d-2i)e_{1},0,ie_1}(d,r)$$

A vector $\alpha$ in $\NN^\infty$ is said to be \textit{odd} if
$(\alpha)_{2i}=0$ for all $i\ge 1$. It follows immediately from the definition of
the multiplicity of a marked floor diagram that  
$C^{\alpha,\beta,\gamma,\delta}(d,r)=0$ if $\alpha$ or $\beta$ is not
odd. 

Given three 
integer numbers
$l\ge 0$, $m\ge 0$ and $r\ge 0$  and four vectors $\alpha$, $\beta$,
$\gamma$ and $\delta$
in $\NN^\infty$, we denote by $\mathcal
S_w(l,m,r,\alpha,\beta,\gamma,\delta)$ the set 
composed by the vectors
$(d_1,\ldots,d_l,k_1,\ldots,k_l,\gamma_1,\ldots,\gamma_l,$
$\delta_1,\ldots,\delta_l,$
$d_1',\ldots,d_m',k_1',\ldots,k_m',$ 
$r_1',\ldots,r_m',$ 
$\alpha_1',\ldots,\alpha_m',\beta_1',
\ldots,\beta_m',$
$\gamma_1',\ldots,\gamma_m',
\delta_1',\ldots,\delta_m')$
in   
$(\NN^*)^{2l}\times (\NN^\infty)^{2l}$
$\times(\NN^*)^{2m}$ $\times \NN^m$
$\times (\NN^\infty)^{4m} $
satisfying 
\begin{itemize}

\item $\forall i, (d_i',k_i',r_i',\alpha_i',\beta_i',\gamma_i',\delta_i')$$\le$$
  (d_{i+1}',k_{i+1}',r_{i+1}',,\alpha_{i+1}',\beta_{i+1}',\gamma_{i+1}',\delta_{i+1}')$ for the lexicographic
order, 

\item $\sum \alpha'_i \le \alpha$,
\item  $\forall i, k'_i$ is odd,
\item $\forall i, \beta'_i\ge e_{k'_i}$, 
\item $\sum (\beta'_i-e_{k'_i})
 =\beta $,
\item $\sum \gamma'_i \le \gamma$,
\item $\sum \delta'_i\le \delta $,

\item $\forall i, I\alpha'_i+I\beta'_i +2I\gamma'_i  +2I\delta'_i=d'_i$.

\item $(d_1,\ldots,d_l,k_1,\ldots,k_l,\gamma_1,\ldots,\gamma_l,$
$\delta_1,\ldots,\delta_l)\in \mathcal S(\frac{d-\sum
  d_i'+1}{2},l,\gamma-\sum \gamma'_i 
  ,\delta-\sum  \delta'_i )$,

\item $\sum (2d_i -1+|\delta_i|) +\sum r'_i = r$.

\item $\forall i, 2d_i'-1+|\beta_i'+2\delta'_i|-r_i'\geq 0$

\end{itemize}

Given two
integer numbers $l\ge 0$ and $r\ge 0$ and four vectors $\alpha$, $\beta$,
$\gamma$ and $\delta$
in $\NN^\infty$, we denote by $\widetilde{\mathcal
S}_w(l,r,\alpha,\beta,\gamma,\delta)$ the set 
composed by the vectors
$(d_1,\ldots,d_l,k_1,\ldots,k_l,\gamma_1,\ldots,\gamma_l,$
$\delta_1,\ldots,\delta_l,$
$d_1',k_1',$ 
$r_1',$ 
$\gamma_1',
\delta_1')$
in   
$(\NN^*)^{2l}\times (\NN^\infty)^{2l}$
$\times(\NN^*)^{2}$$\times \NN $  $\times (\NN^\infty)^{2} $
satisfying 
\begin{itemize}
\item $\gamma'_1 \le \gamma$,
\item $0\le\delta'_1- e_{k'_1} \le\delta $, 
\item $ I\alpha +I\beta +2I\gamma'_1  +2I\delta'_1=d'_1$.

\item $(d_1,\ldots,d_l,k_1,\ldots,k_l,\gamma_1,\ldots,\gamma_l,$
$\delta_1,\ldots,\delta_l)\in \mathcal S(\frac{d-
  d_1'}{2},l,\gamma- \gamma'_1
  ,\delta- \delta'_1 +   e_{k'_1})$,

\item $\sum (2d_i -1+|\delta_i|) + r'_1 = r$.
\item $2d_1'-1+|\beta+2\delta'_1|-r_1'\geq 0$


\end{itemize}

By definition, any element 
$s$     
of $\mathcal S_w(l,m,r,\alpha,\beta,\gamma,\delta)$ 
(resp. $\widetilde{\mathcal S}_w(l,r,\alpha,\beta,\gamma,\delta)$) defines an
element $s'$ of $\mathcal S(\frac{d-\sum
  d_i'+1}{2},l,\gamma-\sum \gamma'_i 
  ,\delta-\sum  \delta'_i )$ (resp. $\mathcal S(\frac{d-
  d_1'}{2},l,\gamma- \gamma'_1
  ,\delta- \delta'_1 +   e_{k'_1})$). We denote by $\sigma(s)$ the
integer $\sigma(s')$ defined in section \ref{ch}.





To any element 
$s$   
of $\mathcal S_w(l,m,r,\alpha,\beta,\gamma,\delta)$,  we associate
the equivalence relation $\simeq_s$ 
on the set $\{1,\ldots ,m\}$ defined by
$$i\simeq_s j \Leftrightarrow (d'_i,  k'_i ,  r'_i,\alpha'_i,\beta'_i,
\gamma'_i, 
 \delta'_i)=(d'_j,  k'_j ,  r'_j,\alpha'_j,\beta'_j,
\gamma'_j, 
 \delta'_j)$$
For each of the equivalent classes of $\simeq_s$, evaluate the factorial of its
cardinal, and denote by $\sigma'(s)$ the product of these factorials.

To any element
$s$   
of $\mathcal S_w(l,m,r,\alpha,\beta,\gamma,\delta)$,
we associate the set $\mathcal
E(s)$ 
 of all $j$ in 
$\{1,\ldots,m\}$ such that $\beta'_j\ge
  e_{k'_j}$  and $ 2d'_j-1+|\beta'_j+2\delta'_j| -2r'_j=1$. 
Given an element 
$j\in\{1,\ldots,m\}$, we denote by $\simeq_s^j$ the restriction to
  $\{1,\ldots,m\}\setminus\{j\}$ of the
equivalence relation $\simeq_s$. For each of the equivalent classes of $\simeq_s^j$, evaluate the factorial of its
cardinal, and denote by $\sigma'(s,j)$ the product of these factorials.

\vspace{2ex}
Let us  introduce 
one more
notation. 
Given an element
$s$     
of $\mathcal S_w(l,m,r,\alpha,\beta,\gamma,\delta)$ or 
$\widetilde{\mathcal S}_w(l,r,\alpha,\beta,\gamma,\delta)$, we
put

$$\begin{array}{lcl}
\Theta(s)&=& \left(\begin{array}{c}
r
\\ 2d_1 -1+|\delta_1|,\ldots, 2d_l -1+|\delta_l|,r'_1,\ldots,r'_m
\end{array} \right) \left(\begin{array}{c}
\alpha
\\ \alpha'_1,\ldots, \alpha'_m
\end{array} \right) 
\left(\begin{array}{c}
\gamma
\\  \gamma_1,\ldots,\gamma_l,\gamma'_1,\ldots, \gamma'_m
\end{array} \right)
\\ 
\\&& \prod_{i=1}^l (-1)^{d_i}(\delta_i)_{k_i}k_i2^{2d_i
  -2+|\delta_i+\gamma_i|} N^{\gamma_i,\delta_i}(d_i)
\end{array}$$
where if $s$ is in  $\widetilde{\mathcal
  S}_w(l,r,\alpha,\beta,\gamma,\delta)$, then there are no $\alpha_i'$ 
elements and we set 
the value of
the corresponding multinomial coefficient equal to 1.

\begin{thm}\label{form}
The numbers $C^{\alpha,\beta,\gamma,\delta}(d,r)$, with $\alpha$ and
$\beta$ odd, are given by the initial values
$C^{e_1,0,0,0}(1,0)=C^{0,e_1,0,0}(1,1)=1$ and the relations

\vspace{2ex}
(1) if $2d-1 +|\beta +2\delta | -2r >0 $ 
\vspace{2ex}

$\begin{array}{lcl}
C^{\alpha, \beta,\gamma ,\delta}(d,r) &=&\sum_{k\ odd|\beta\ge e_k} C^{\alpha+e_k,
 \beta-e_k, \gamma ,\delta}(d,r) \ \ + 
\end{array}$

$\begin{array}{lcl}
\\\hspace{14ex} && \sum_{\begin{array}{ll}_{l,m\ge 0}\\ _{s\in\mathcal
S_w(l,m,r,\alpha,\beta,\gamma,\delta)}\end{array}}
\left[\frac{\Theta(s)}{\sigma(s)\sigma'(s)} \prod_{i=1}^m
(\beta'_i)_{k'_i}C^{\alpha'_i,\beta'_i,\gamma_i',\delta'_i}(d'_i,r'_i)
\right.

\\&&\left. \left(\begin{array}{c}
2d-2 +|\beta +2\delta | -2r
\\ 2d'_1 -1+|\beta'_1+ 2\delta'_1|-2r'_1,\ldots, 2d'_m -1+|\beta'_m+ 2\delta'_m|-2r'_m
\end{array} \right)
\right] 
\end{array}
$

\vspace{6ex}

(2) if $2d-1 +|\beta +2\delta | -2r =0 $ 
\vspace{2ex}

$\begin{array}{lcl}
C^{\alpha, \beta,\gamma ,\delta}(d,r) &=&\sum_{k|\delta\ge e_k}kC^{\alpha,
 \beta, \gamma+e_k ,\delta-e_k}(d,r-1) \ \ +
\\
\\ &&\sum_{\begin{array}{l}_{l,m\ge0}\\ _{K\textrm{ odd}| \beta\ge e_K}\\
_{s\in \mathcal S_w(l,m,r-1,\alpha,\beta-e_K,\gamma,\delta)}\end{array}}
\frac{K\Theta(s)}{\sigma(s)\sigma'(s)} \prod_{i=1}^m
(\beta'_i)_{k'_i}C^{\alpha'_i,\beta'_i,\gamma'_i,\delta'_i}(d'_i,r'_i)

\\ &+& \sum_{\begin{array}{l} _{l,m\ge 0}\\ _{s\in \mathcal
S_w(l,m,r-1,\alpha,\beta,\gamma,\delta)} \\ _{j \in \mathcal E(s)}\end{array}}
\frac{\Theta(s)k'_jC^{\alpha'_j+e_{ k'_j},\beta'_j-e_{ k'_j},\gamma'_j,\delta'_j}(d'_j,r'_j)}{\sigma(s)\widetilde{\sigma}'(s,j)}

 \prod_{i=1,i\ne j}^m (\beta'_i)_{k'_i}C^{\alpha'_i,\beta'_i,\gamma'_i,\delta'_i}(d'_i,r'_i)


\\ &-& \sum_{\begin{array}{l}_{l\ge 0}\\ _{s\in\widetilde{\mathcal
S}_w(l,r-1,\alpha,\beta,\gamma,\delta)}\end{array}}
\frac{\Theta(s)}{\sigma(s)} 2^{|\gamma  - \gamma'_1 -\sum_{i=1}^l\gamma_i|+1}(\delta'_1)_{k'_1}
k'_1C^{\alpha,\beta,\gamma'_1,\delta'_1}(d'_1,r'_1)
\end{array}
$

\end{thm}

\begin{rem}
If one plugs $r=0$ in Theorem \ref{form}, then one finds 
the
Itenberg-Kharlamov-Shustin formula from \cite{IKS3}.
\end{rem}

\section{Proof of Theorem \ref{form}}\label{proof}

The proof of Theorem \ref{form} 
is
divided into two parts,
depending on wether $2d-1+|\beta+2\delta|-2r>0$ or not. 
In this whole section
$\alpha$ and $\beta$ are two odd vectors in $\NN^\infty$.

\subsection{The case $2d-1+|\beta+2\delta|-2r>0$}


By definition, one has
$$C^{\alpha,\beta,\gamma,\delta}(d,r)=\sum_{(\mathcal D,\m)\in
   A}\mu_r^\RR(\mathcal D,\m) + \sum_{(\mathcal D,\m)\in
   B}\mu_r^\RR(\mathcal D,\m)$$
where $A$ (resp. $B$) 
consists
of all $r$-real
marked floor diagrams of degree $d$ 
 and type $(\alpha,\beta,\gamma,\delta)$ such that
 $\m(\left|\alpha+2\gamma\right|+1)\in\text{Edge}^{\infty}(\mathcal
 D)$ (resp. $\m(\left|\alpha+2\gamma\right|+1)\in\text{Vert}(\mathcal
 D)$). 

\subsubsection{Marked floor diagrams in $A$}

Let us de denote by $A_k $ the set of  all $r$-real marked floor
diagrams of degree $d$ 
 and type
 $(\alpha+e_k,\beta-e_k,\gamma,\delta)$.
 There exists a bijection $\Phi$ from the set $A$ to the union of
all the sets $A_k$
for  $k$  such that
 $\beta\ge e_k$.  If  $(\mathcal D,\m)$ is a marked floor diagram in $A$ and
$k$ is the weight of the edge $\m(\left|\alpha+2\gamma\right|+1)$,
we define $\Phi(\mathcal{D},\m)= (\mathcal{D},\m')$ where 

\begin{itemize} 
\item $\m'(i)=\m(i)$ if $i\leq \sum_{j=1}^{k}\alpha_j$ or $i\geq \left|\alpha+2\gamma\right|+2$,
\item $\m'(\sum_{j=1}^{k}\alpha_j+1)$ is the source adjacent to $\m(\left|\alpha+2\gamma\right|+1)$,
\item $\m'(i)=\m(i-1)$ if $\sum_{j=1}^{k}\alpha_j+2\leq i\leq \left|\alpha+2\gamma\right|+1$.
\end{itemize}

By definition of the multiplicity, one has $\mu^\RR_r(\mathcal 
D,\m)=\mu^\RR_r(\mathcal D,\m')$, and then
 $$\sum_{(\mathcal D,\m)\in A} \mu^\RR_r(\mathcal D
,\m)=\sum_{k|\beta\ge e_k}
\sum_{(\mathcal D,\m')\in A_k} \mu^\RR_r(\mathcal D,\m')$$ 

Note that
$k$ is odd since $\beta$ is odd, hence one has

$$\sum_{(\mathcal D,\m)\in A} \mu^\RR_r(\mathcal D ,\m)=\sum_{k\ odd
  |\beta\ge e_k}C^{\alpha+e_k,\beta-e_k,\gamma,\delta}(d,r)$$

\subsubsection{Marked floor diagrams in $B$}\label{BB}
Let us first explain backward the idea to obtain the second sum in
the right-hand side 
of the formula in Theorem
\ref{form},  equation (1).
If  $(\mathcal D,\m)$ is a marked floor diagram in $B$, by ``cutting'' the vertex
 $\m(\left|\alpha+2\gamma\right|+1)$ from $\mathcal D$, one obtains
 several marked floor diagrams of genus $0$ and of lower
 degrees. Since $(\mathcal D,\m)$ is  $r$-real,  the new marked
 floor diagrams  contained in $\Im(\mathcal D,\m,r)$ are naturally
 coupled in pairs by the 
 involution $\rho_{\mathcal D,\m,r}$ (see definition
 \ref{defi real}). Moreover,  none of
the edges in $\text{Edge}^{\infty}(\mathcal 
 D)$ adjacent to the vertex $\m(\left|\alpha+2\gamma\right|+1)$ is in
 the image of $\m$.

\vspace{2ex}
Let $l$ and $m$ be two nonnegative integer numbers such that the set
$\mathcal S_w(l,m,r,\alpha,\beta,\gamma,\delta)$ is not empty. For $s$
in $\mathcal 
S_w(l,m,r,\alpha,\beta,\gamma,\delta)$, denote by $\mathfrak P(s)$ the set of
all  $2l+m$-tuples of marked floor diagrams 
$(\mathcal D_1,\m_1),\dots,(\mathcal D_{2l},\m_{2l}),(\mathcal
D'_1,\m'_1),\dots,(\mathcal D'_m,\m'_m)$
 such that 

\begin{itemize}
\item   $(\mathcal D_{2i-1},\m_{2i-1})$ and  $(\mathcal
  D_{2i},\m_{2i})$  are two equivalent   marked floor diagrams  of
  degree $d_i$
  and type 
  $(\gamma_i,\delta_i,0,0)$,  
\item $(\mathcal D'_i,\m'_i)$ is $r'_i$-real of degree $d'_i$
and type $(\alpha'_i,\beta'_i,\gamma'_i,\delta'_i)$.
\end{itemize}
Let  $\phi_i$ be a homeomorphism of 
the oriented graph identifying 
$(\mathcal D_{2i-1},\m_{2i-1})$ and  $(\mathcal
  D_{2i},\m_{2i})$.

Starting from an element of $ \mathfrak P(s)$, we construct several
elements of $B$ in the following way 

\begin{enumerate}

\item For all $i$ in $\{1,\ldots, l\}$ choose an element $a_i$ of
$\text{Edge}^{\infty}(\mathcal D_{2i-1})$ which is in the image of
  $\m_{2i-1}$ and of weight
$k_i$.
 Since
$\delta_i\ge e_{k_i}$, it is always
possible to choose such an $a_i$.

\item
For all $i$ in $\{1,\dots, m\}$ choose an element $a_i'$ of
$\text{Edges}^{\infty}(\mathcal 
D'_i)$ which is in the image of $\m'_i$ but not in $\Im(\mathcal
D'_i,\m'_i,r'_i)$, and of weight $k'_i$. Since
$\beta'_i\ge e_{k'_i}$, it is always
possible to choose such an $a_i'$.

\item Construct a new oriented tree  $\widetilde{\mathcal D}$ out of
  $(\mathcal D_1,\m_1),\dots,(\mathcal D_{2l},\m_{2l}),(\mathcal
D'_1,\m'_1),\dots,(\mathcal D'_m,\m'_m)$ by
identifying all the sources adjacent to the edges $a_i$, $\phi_i(a_i)$,
and $a_j'$. Denote
  this vertex by $v$.

\item By adding sources and edges adjacent to the vertex $v$, complete
  $\widetilde 
  {\mathcal D}$ into a (unique) floor diagram 
  $\mathcal D$ of  degree $d$, genus $0$, with $(\alpha)_j +
  (\beta)_j + 2(\gamma)_j + 2(\delta)_j$ edges in
  $\text{Edges}^{\infty}(\mathcal  
D)$ of weight $j$ for all $j\ge 1$. Denote by $v_1,\ldots, v_t$ the
sources added.

\item Define $\alpha'_{m+1}=\alpha -\sum_{i=1}^m\alpha'_i$ and
 $\gamma'_{m+1}=\gamma - \sum_{i=1}^l\gamma_i- \sum_{i=1}^m\gamma'_i$.

 \item For all $j\ge 1$, choose a partition $(I_i^j)_{1\leq i\leq m+1}$
  of the set  $ \{1,\dots,(\alpha)_j\}$ 
   such that $\sharp I_i^j=(\alpha'_i)_j$ for
  all $i$.

\item  For all $j\ge 1$, choose a partition $(\widehat
  I_i^j)_{1\leq i\leq l}\cup (\widetilde
  I_i^j)_{1\leq i\leq m+1}$
  of the set  $ \{1,\dots,(\gamma)_j\}$ 
   such that $\sharp \widehat I_i^j=(\gamma_i)_j$ and $\sharp
   \widetilde I_i^j=(\gamma'_i)_j$ for
  all $i$.

\item Choose a partition $(J_i)_{1\leq i\leq
  m}$ of the set $\{1,\dots,
  2d-2+\left|\beta+2\delta\right|-2r\}$ such that 
$\sharp J_i=2d'_i-1+\left|\beta'_i+2\delta'_i\right|-2r'_i$ for
  all $i$.

\item Choose a partition $(\widehat J_i)_{1\leq i\leq
  l}\cup (\widetilde J_i)_{1\leq i\leq
  m}$ of the set $\{1,\dots,
  r\}$ such that 
$\sharp \widehat J_i=2d_i-1+|\delta_i|$ and
 $\sharp \widetilde J_i=r'_i$ for
 all $i$.

\item For all $i$ in $\{1,\ldots, l\}$, choose
  a vector $\varepsilon_i$ 
  in $\{0,1\}^{2d_i-1+|\gamma_i+\delta_i|}$.

\item Choose a marking $\m$ of $\mathcal D$ of type
  $(\alpha,\beta,\gamma,\delta)$ such that  
\begin{enumerate}

\item $\m(|\alpha|+ 2|\gamma|+1)=v$,

\item for all  $j\ge 1$ and all $k$ in $I_{m+1}^j$, 
$\m(\sum_{t=1}^{j-1}\left(\alpha)_t+k\right)$ is a 
  source $v_q$ (see step (4)) of $\mathcal D$ of  divergence $-j$
  (note that different 
  choices of $v_q$ produce equivalent marked floor diagrams),

\item for all $j\ge 1$ and all $k$ in $\widetilde I_{m+1}^j$, 
  $\m\left(|\alpha|+2\sum_{t=1}^{j-1}(\gamma)_t+2k-1\right)$ and
  $\m\left(|\alpha|+2\sum_{t=1}^{j-1}(\gamma)_t\right.$ $+2k\Big)$ are two sources $v_q$
  and $v_{q'}$  of
  $\mathcal D$ of  divergence $-j$  (again, different 
  choices of $v_q$ and $v_{q'}$ produce equivalent marked floor diagrams), 

\item for all $j\ge 1$ and all $i$ in $\{1,\ldots,m\}$, if $k$ is the $h$-th
  element (for the natural ordering on $I_i^j$) of 
  $I_j^i$, then
$$\m\left (\sum_{t=1}^{j-1}(\alpha)_t+k\right)=
  \m'_j\left(\sum_{t=1}^{j-1}(\alpha'_j)_t+h\right)$$    

\item for all $j\ge 1$ and all $i$ in $\{1,\ldots,l\}$, if $k$ is the
  $h$-th element of $\widehat
  I_i^j$, then
 $$\m\left(|\alpha|+2\sum_{t=1}^{j-1}(\gamma)_t+2k-1+(\varepsilon_i)_{\sum_{t=1}^{j-1}(\gamma_i)_t
    + h}
  \right)=
\m_{2i-1}\left(\sum_{t=1}^{j-1}(\gamma_i)_t+h \right) $$

and
$$
  \m\left(|\alpha|+2\sum_{t=1}^{j-1}(\gamma)_t+2k-
  (\varepsilon_i)_{\sum_{t=1}^{j-1}(\gamma_i)_t 
    + h}
  \right)=
  \phi_i\circ \m_{2i-1}\left(\sum_{t=1}^{j-1}(\gamma_i)_t+h\right)$$ 

\item  for all $j\ge 1$ and all $i$ in $\{1,\ldots,m\}$, if $k$ is the
  $h$-th element of $\widetilde
  I_i^j$, then

 $$\m\left(|\alpha|+2\sum_{t=1}^{j-1}
  (\gamma)_t+2k-1\right)=\m'_i\left(|\alpha'_i|+2\sum_{t=1}^{j-1} 
  (\gamma'_i)_t+2h-1\right)$$ 
  and 
$$\m\left(|\alpha|+2\sum_{t=1}^{j-1}(\gamma)_t+2k\right)=
\m'_i\left(|\alpha'_i|+2\sum_{t=1}^{j-1}(\gamma'_i)_t+2h\right)$$  

\item for all $i$ in  $\{1,\ldots,m\}$, if $k$ is the $h$-th element
  of $J_i$, then
$$\m\left(|\alpha|+ 
  2|\gamma|+k+1\right)=\m'_i\left(|\alpha'_i|+ 2|\gamma'_i|+h\right)$$ 

\item for all $i$ in  $\{1,\ldots,l\}$, if $k$ is the $h$-th element
  of $\widehat J_i$, then (recall that $n=2d-1+|\alpha+\beta+2\gamma+2\delta|$)
$$\m\left(n-2k+1 + (\varepsilon_i)_{|\gamma_i|+
    h}\right)=\m_{2i-1}\left(2d_i-1+|\delta_i|-h+1\right) $$ 
and 
$$\m\left(n-2k+2- (\varepsilon_i)_{|\gamma_i|+
    h}\right)=\phi_i\circ \m_{2i-1}\left(2d_i-1+|\delta_i|-h+1\right) $$ 

\item for all $i$ in  $\{1,\ldots,m\}$, if $k$ is the $h$-th element of
  $\widetilde J_i$, then
 $$\m\left(n-2k+1\right)=\m'_i\left(2d'_i-1+|\beta'_i+2\delta'_i|-2h+1\right)$$

  and 

 $$\m\left(n-2k+2\right)=\m'_i\left(2d'_i-1+|\beta'_i+2\delta'_i|-2h+2\right)$$

\end{enumerate}

\end{enumerate}

In this way, we construct all marked floor diagrams in $B$. Moreover,
two distinct elements of $\mathcal S_w(l,m,r,\alpha,\beta,\gamma,\delta)$ will
produce non equivalent marked 
floor diagrams. If a marked floor diagram  $(\mathcal D ,\m)$ of type
  $(\alpha,\beta,\gamma,\delta)$ is constructed out of an element $s$
of $\mathcal S_w(l,m,r,\alpha,\beta,\gamma,\delta)$, then $(\mathcal D ,\m)$ is
obtained exactly $2^l\sigma(s)\sigma'(s)$ 
times. Indeed, since $(\mathcal D_{2i-1} ,\m_{2i-1})$ and $(\mathcal
D_{2i} ,\m_{2i})$ are equivalent, the two vectors $\varepsilon_i$
and $(1,\ldots,1)- \varepsilon_i$ produce 2 equivalent marked floor diagrams.
The factor $\sigma(s)\sigma'(s)$ is clear.

By construction, one has

$$\mu^\RR_r(\mathcal D ,\m)=\prod_{i=1}^l(-1)^{d_i}k_i\mu^\CC(\mathcal
D_{2i} ,\m_{2i})\prod_{i=1}^m\mu^\RR_{r'_i}(\mathcal D'_i ,\m'_i)$$ 

Now,  the second sum in
the right-hand side 
of the formula in Theorem
\ref{form},  equation (1), follows from 
all possible choices in our
construction,
and from Theorem \ref{NFD}.

\subsection{The case $2d-1+|\beta+2\delta|-2r=0$}

 In this case, if $(\mathcal D,\m)$ is an $r$-real marked floor
diagram of degree $d$
and type
$(\alpha,\beta,\gamma,\delta)$, then 
$\{\left|\alpha+2\gamma\right|+1,\left|\alpha+2\gamma\right|+2\}$
is an $r$-pair of $(\mathcal D,\m)$. The four terms in the
right-hand side of 
the formula in Theorem
\ref{form},  equation (2), come from  
consideration of four subcases.
Let $A'$, $B'$, $C'$ and $D'$ be the 
the sets of all $r$-real marked floor
diagrams of degree $d$
and type
$(\alpha,\beta,\gamma,\delta)$  satisfying respectively 

\begin{itemize}
\item  both
$\m\left(\left|\alpha+2\gamma\right|+1\right)$ and
$\m\left(\left|\alpha+2\gamma\right|+2\right)$ are in
$\text{Edge}^{\infty}(\Gamma)$,   
\item 
$\m\left(\left|\alpha+2\gamma\right|+1\right)$ is in
$\text{Edge}^{\infty}(\Gamma)$  and
$\m\left(\left|\alpha+2\gamma\right|+2\right)$ is in $\text{Vert}(\Gamma)$,
\item
$\m\left(\left|\alpha+2\gamma\right|+1\right)$ is in $ \text{Vert}(\Gamma)$ and
$\m\left(\left|\alpha+2\gamma\right|+2\right)$ is in
$\text{Edge}(\Gamma)$, 
\item both
$\m\left(\left|\alpha+2\gamma\right|+1\right)$ and
$\m\left(\left|\alpha+2\gamma\right|+2\right)$ are in $\text{Vert}(\Gamma)$.
\end{itemize}

\subsubsection{Marked floor diagrams in $A'$}

There exists a bijection $\Phi'$ from the set $A'$ to the union of
 all $(r-1)$-real marked floor diagrams of degree $d$
 and type
 $(\alpha,\beta,\gamma+e_k,\delta-e_k)$ for  $k$  such that
 $\delta\ge e_k$.  If  $(\mathcal D,\m)$ is a marked floor diagram in $A'$ and
$k$ is the weight of the edge $\m(\left|\alpha+2\gamma\right|+1)$,
we define $\Phi'(\mathcal{D},\m)= (\mathcal{D},\m')$ where 

\begin{itemize} 
\item $\m'(i)=\m(i)$ if $i\leq |\alpha|+2\sum_{j=1}^{k}\gamma_j$ or $i\geq \left|\alpha+2\gamma\right|+3$,
\item $\m'(|\alpha|+2\sum_{j=1}^{k}\gamma_j+1)$ is the source adjacent to $\m(\left|\alpha+2\gamma\right|+1)$,
\item $\m'(|\alpha|+2\sum_{j=1}^{k}\gamma_j+2)$ is the source adjacent to $\m(\left|\alpha+2\gamma\right|+2)$,
\item $\m'(i)=\m(i-2)$ if $|\alpha|+2\sum_{j=1}^{k}\gamma_j+3\leq
  i\leq \left|\alpha+2\gamma\right|+2$. 
\end{itemize}

Note that since $(\mathcal D,\m)$ is $r$-real, the weight of the
edges $\m(\left|\alpha+2\gamma\right|+1)$ and
$\m(\left|\alpha+2\gamma\right|+2)$ is
 the same. 
One has $\mu^\RR_r(\mathcal
D,\m)=k\mu^\RR_r(\mathcal D,\m')$, hence the first sum of the
right-hand side of 
the formula in Theorem
\ref{form},  equation (2), is given by $\sum_{(\mathcal D,\m)\in A'}
\mu^\RR_r(\mathcal D ,\m)$.

\subsubsection{Marked floor diagrams in $B'$}

Choose $K\ge 1$ such that $\beta\ge e_K$. 
Let $l$ and $m$ be two nonnegative integers such that the set
$S_w(l,m,d,k,r-1,\alpha,\beta-e_K,\gamma,\delta)$ is not empty, and define
the set $\mathfrak P(s)$ as in section \ref{BB} for any $s$ in
$S_w(l,m,d,k,r-1,\alpha,\beta-e_K,\gamma,\delta)$. 
 Then, starting from an element of $\mathfrak P(s)$ we construct several
 elements of $B'$ as in the step (1) - (11) of section \ref{BB},
 except for the following modifications

\begin{itemize}
\item[(8B')] Define $J_i=\emptyset$ for all $i$ in $\{1,\ldots,m\}$.
\item[(9B')] Choose a partition $(\widehat J_i)_{1\leq i\leq
  l}\cup (\widetilde J_i)_{1\leq i\leq
  m}$ of the set $\{1,\dots,
  r-1\}$ such that 
$\sharp \widehat J_i=2d_i-1+|\delta_i|$ and
 $\sharp \widetilde J_i=r'_i$ for
 all $i$.
\item[(11B')] (a) 
 $\m(|\alpha|+ 2|\gamma|+2)=v$, and $\m(|\alpha|+
  2|\gamma|+1)$ is an edge in $\text{Edge}^\infty(\mathcal D)$ of
  weight $K$ adjacent to $v$

\end{itemize}

By construction, one has

$$\mu^\RR_r(\mathcal D ,\m)=K \prod_{i=1}^l(-1)^{d_i}k_i\mu^\CC(\mathcal
D_{2i} ,\m_{2i})\prod_{i=1}^m\mu^\RR_{r'_i}(\mathcal D'_i ,\m'_i)$$ 

Now,  the second sum in
the right-hand side 
of the formula in Theorem
\ref{form},  equation (2), follows from all possible choices in our
construction.


\subsubsection{Marked floor diagrams in $C'$}

Let $l$ and $m$ be two nonnegative integers such that the set
$S_w(l,m,d,k,r-1,\alpha,\beta,\gamma,\delta)$ is not empty, and define
the set $\mathfrak P(s)$ as in section \ref{BB} for any $s$ in
$S_w(l,m,d,k,r-1,\alpha,\beta,\gamma,\delta)$. 
 Then, starting from an element of $\mathfrak P(s)$ we construct several
 elements of $C'$ as in the step (1) - (11) of section \ref{BB},
 except for the following modifications

\begin{itemize}
\item[(0C')] Choose $j$ in $\{1,\ldots, m\}$ such that $\beta_j\ge k'_j$, the
  weight of $\m'_j(|\alpha'_j|+ 2|\gamma'_j|+1)$ is $k'_j$, and
  $2d'_j-1 +|\beta'_j+2\gamma'_j|-2r'_j=1$. 
\item[(2C')] For all $1\leq i\leq m$, $i\neq j$, choose an element $a_i'$ of
$\text{Edges}^{\infty}(\mathcal 
D'_i)$ which is in the image of $\m'_i$ but not in $\Im(\mathcal
D'_i,\m'_i,r'_i)$, and of weight $k'_i$. 
\item[(8C')] Define $J_i=\emptyset$ for all $i$ in $\{1,\ldots,m\}$.
\item[(9C')] Choose a partition $(\widehat J_i)_{1\leq i\leq
  l}\cup (\widetilde J_i)_{1\leq i\leq
  m}$ of the set $\{1,\dots,
  r-1\}$ such that 
$\sharp \widehat J_i=2d_i-1+|\delta_i|$ and
 $\sharp \widetilde J_i=r'_i$ for
 all $i$.
\item[(11C')] (a) 
 $\m(|\alpha|+ 2|\gamma|+1)=v$, and $\m(|\alpha|+
  2|\gamma|+2)=\m'_j(|\alpha'_j|+ 2|\gamma'_j|+1)$.
  
\end{itemize}
By construction, one has

$$\mu^\RR_r(\mathcal D ,\m)=k'_j\prod_{i=1}^l(-1)^{d_i}k_i\mu^\CC(\mathcal
D_{2i} ,\m_{2i})\prod_{i=1}^m\mu^\RR_{r'_i}(\mathcal D'_i ,\m'_i)$$ 
 
Now,  the third sum in
the right-hand side 
of the formula in Theorem
\ref{form},  equation (2), follows from all possible choices in our
construction.   

\subsubsection{Marked floor diagrams in $D'$}

In this case, by "cutting" the vertices
 $\m(\left|\alpha+2\gamma\right|+1)$ and
$\m(\left|\alpha+2\gamma\right|+2)$ from $\mathcal D$, one obtains 
 several marked floor diagrams of genus $0$ and of lower 
 degrees. Since $\mathcal D$ is a tree and is $r$-real, exactly one of
 these marked floor 
 diagrams is adjacent to both cut vertices, and all the other  are naturally 
 coupled in pairs by the map $\rho_{\mathcal D,\m,r}$ (see definition
 \ref{defi real}). 
Moreover,  any edge in $\text{Edge}^{\infty}(\mathcal 
 D)$ adjacent to the vertex $\m(\left|\alpha+2\gamma\right|+1)$ is not 
in 
 the image of $\m$ and is naturally coupled  to an edge
adjacent to the vertex $\m(\left|\alpha+2\gamma\right|+2)$ since the
diagram is $r$-real.
In particular, both edges have
the same weight.
 
\vspace{2ex}
Let $l$ be a nonnegative integer such that the set
$\widetilde S_w(l,r-1,\alpha,\beta,\gamma,\delta)$ is not empty.
For $s$
in this set, 
  define
the set $\mathfrak P(s)$ as in section \ref{BB} with $m=1$.

Starting from an element of $ \mathfrak P(s)$, we construct several
elements of $D'$ in the following way 

\begin{enumerate}

\item For all $0\leq i\leq l$ choose an element $a_i$ of
$\text{Edge}^{\infty}(\mathcal D_{2i-1})$  which is in the image of $\m_{2i-1}$
  and of weight
$k_i$.

\item Choose an edge $a'_1$ of
$\mathcal D'_{1}$ in $\text{Edge}^{\infty}(\mathcal D'_1)\cap \Im(\mathcal D'_{1},m'_1,r'_1)$ of weight
$k'_1$.

\item Construct a new oriented tree  $\widetilde{\mathcal D}$ out of
  $(\mathcal D_1,\m_1),\dots,(\mathcal D_{2l},\m_{2l}),(\mathcal
D'_1,\m'_1)$ by
identifying  
\begin{enumerate}
\item  all the sources adjacent to the edges
  $a_i$ and $a'_1$

\item  all the
  sources adjacent to the edges 
  $\phi_i(a_i)$ and $\rho_{\mathcal D'_{1},m'_1,r'_1}(a'_1)$.
\end{enumerate}
Denote by $v$ and
 $v'$ the 2 vertices added. 

\item Construct a degree $d$ and genus
  $0$ floor diagram  $\mathcal D$ out of $\widetilde{\mathcal D}$ by
adding sources $v_1,\dots, v_t, v'_1,$ $\dots,$ $v'_t$ and edges 
$(v_1,v),\dots ,$ $ (v_t,v),$ $ (v'_1,v'),$ $ \dots,$ $(v'_t,v')$, such
that $\mathcal D$ has $(\alpha)_j +
   (\beta)_j + 2(\gamma)_j + 2(\delta)_j$ edges in
   $\text{Edges}^{\infty}(\mathcal  
 D)$ of weight $j$ for all $j\ge 1$, and such that there are as
 many edges $(v_i,v)$ of weight $j$  
 as edges $(v'_i,v')$ of weight $j$
 for all $j\ge 1$. 


\item Define $\gamma'_{2}=\gamma - \sum_{i=1}^l\gamma_i- \gamma'_1$.

\item  For all $j\ge 1$, choose a partition $(\widehat
  I_i^j)_{1\leq i\leq l}\cup (\widetilde
  I_1^j)\cup (\widetilde
  I_2^j)$
  of the set  $ \{1,\dots,(\gamma)_j\}$ 
   such that $\sharp \widehat I_i^j=(\gamma_i)_j$ and $\sharp
   \widetilde I_i^j=(\gamma'_i)_j$.

\item Choose a partition $(\widehat J_i)_{1\leq i\leq
  l}\cup (\widetilde J_1)$ of the set $\{1,\dots,
  r-1\}$ such that 
$\sharp \widehat J_i=2d_i-1+|\delta_i|$ and
 $\sharp \widetilde J_1=r'_1$.

\item Choose
  a vector $\varepsilon$ 
  in $\{0,1\}^{|\gamma'_2|}$.
  
\item For all $i$ in $\{1,\ldots, l\}$, choose
  a vector $\varepsilon_i$ 
  in $\{0,1\}^{2d_i-1+|\gamma_i+\delta_i|}$.

\item Choose a marking $\m$ of $\mathcal D$ of type
  $(\alpha,\beta,\gamma,\delta)$ such that  
\begin{enumerate}

\item $\m(|\alpha|+ 2|\gamma|+1)=v$ and $\m(|\alpha|+ 2|\gamma|+2)=v'$,

\item for all $j\ge 1$ and if $k$ is the $h$-th element of $\widetilde
  I_{2}^j$, then  
  $\m\left(|\alpha|+2\sum_{t=1}^{j-1}(\gamma)_t+2k-1 \right.$
  $\left.+(\varepsilon)_{\sum_{t=1}^{j-1}(\gamma'_i)_t+h}\right)$ and 
  $\m\left(|\alpha|+2\sum_{t=1}^{j-1}(\gamma)_t+2k \right.$
  $\left.-(\varepsilon)_{\sum_{t=1}^{j-1}(\gamma'_i)_t+h}\right)$
  are respectively a source $v_q$ and $v'_{q'}$ of divergence $-j$,  

\item If $1\leq i\leq |\alpha|$, 
$\m(i)=\m'_1(i)$,

\item for all $j\ge 1$ and all $i$ in $\{1,\ldots,l\}$, if $k$ is the
  $h$-th element of $\widehat
  I_i^j$, then
 $$\m\left(|\alpha|+2\sum_{t=1}^{j-1}(\gamma)_t+2k-1+(\varepsilon_i)_{\sum_{t=1}^{j-1}(\gamma_i)_t
    + h}
  \right)=
\m_{2i-1}\left(\sum_{t=1}^{j-1}(\gamma_i)_t+h \right) $$

and
$$
  \m\left(|\alpha|+2\sum_{t=1}^{j-1}(\gamma)_t+2k-
  (\varepsilon_i)_{\sum_{t=1}^{j-1}(\gamma_i)_t 
    + h}
  \right)=
  \phi_i\circ \m_{2i-1}\left(\sum_{t=1}^{j-1}(\gamma_i)_t+h\right)$$ 

\item  for all $j\ge 1$ , if $k$ is the
  $h$-th element of $\widetilde
  I_1^j$, then

 $$\m\left(|\alpha|+2\sum_{t=1}^{j-1}
  (\gamma)_t+2k-1\right)=\m'_1\left(|\alpha'_1|+2\sum_{t=1}^{j-1} 
  (\gamma'_1)_t+2h-1\right)$$ 
  and 
$$\m\left(|\alpha|+2\sum_{t=1}^{j-1}(\gamma)_t+2k\right)=
\m'_1\left(|\alpha'_1|+2\sum_{t=1}^{j-1}(\gamma'_1)_t+2h\right)$$

\item for all $i$ in  $\{1,\ldots,l\}$, if $k$ is the $h$-th element
  of $\widehat J_i$, then 
$$\m\left(n-2k+1 + (\varepsilon_i)_{|\gamma_i|+
    h}\right)=\m_{2i-1}\left(2d_i-1+|\delta_i|-h+1\right) $$ 
and 
$$\m\left(n-2k+2- (\varepsilon_i)_{|\gamma_i|+
    h}\right)=\phi_i\circ \m_{2i-1}\left(2d_i-1+|\delta_i|-h+1\right) $$

\item if if $k$ is the $h$-th element
  of $\widetilde J_1$,
  
  $$\m\left(n-2k+1\right)=\m'_1\left(2d'_1-1+|\beta'_1+2\delta_i|-2h+1\right) $$ 
and 
$$\m\left(n-2k+2\right)=\m'_1\left(2d'_1-1+|\beta'_1+2\delta_i|-2h+2\right) $$ 

\end{enumerate}

\end{enumerate}

In this way, we construct all marked floor diagrams in $D'$. Moreover,
any element of $D'$  is
obtained exactly $\sigma(s)\sigma'(s)$ 
times for some $s$ in $\widetilde S_w(l,r-1,\alpha,\beta,\gamma,\delta)$.

By construction, one has

$$\mu^\RR_r(\mathcal D ,\m)=-k'_1\mu^\RR_{r'_1}(\mathcal D'_1
,\m'_1)\prod_{i=1}^l(-1)^{d_i}k_i\mu^\CC(\mathcal 
D_{2i} ,\m_{2i})$$ 

Now,  the fourth sum in
the right-hand side 
of the formula in Theorem
\ref{form},  equation (2), follows from all possible choices in our
construction.

\section{Further remarks}\label{remark}

\subsection{More recursive formulas}

Using floor diagrams, one can enumerate
complex curves in other complex  varieties than $\CC P^2$, for example
in Hirzebruch surfaces, or in $\CC P^n$ blown up in a small number of
points.
Floor diagrams 
 allow one to  compute Welschinger invariants  
as soon as 
they are defined
and a theorem similar to Theorem \ref{WFD} holds
 (see \cite{Br6b} and \cite{Br7}).
Hence, 
 one can easily write  formulas 
in the style of
 Theorem \ref{form}
for other varieties than $\CC P^2$, using the same technic as here.
As an example, we give below a
  formula to compute Weslchinger invariants of $\RR P^3$ for rational
  curves passing through a configuration of real points. This
  invariants 
are
defined in \cite{Wel2}, and 
were
computed for the
  first time in \cite{Br6} and \cite{Br7} using floor diagrams.

Let us denote by $W_3(d)$ the Welschinger invariant of degree $d$
of $\RR P^3$ for configurations of $2d$ real points. Given an
integer number $l\ge 0$ and two odd vectors $\alpha$ and $\beta$ in
$\NN^\infty$, we denote by 
$\mathcal S_3(l,\alpha,\beta)$ the set 
composed by the vectors
$(d_1,\ldots,d_l,$
$k_1,\ldots,k_l,$
$\alpha_1,\ldots,\alpha_l,$
$\beta_1,
\ldots,\beta_l)\in   
(\NN^*)^{2l}\times (\NN^\infty)^{2l}$
satisfying 

\begin{itemize}
\item $\forall i, (d_i,k_i,\alpha_i,\beta_i)\le
  (d_{i+1},k_{i+1},\alpha_{i+1},\beta_{i+1})$ for the lexicographic order, 
\item $\sum d_i < I\alpha + I\beta$,
\item $\sum  \alpha_i \le \alpha$,
\item $\forall i, \ k_i$ is odd,
\item $\forall i, \beta_i\ge e_{k_i}$, 
\item $\sum  (\beta_i-e_{k_i})
 =\beta $,

\item $\forall i, I\alpha_i +I\beta_i=d_i$,

\item $l+|\alpha -\sum \alpha_i|=3d - 3\sum d_i-2.$

\end{itemize}

Given an integer number $d\ge 1$ and  two odd
vectors $\alpha$ and $\beta$in $\NN^\infty$ satisfying $I\alpha
+I\beta=d$, we define the 
numbers $ W_3^{\alpha,\beta}(d)$ by
the initial value 
$W_3^{e_1,0}(1)=1$ and the relation

 $$\begin{array}{lcl}
W_3^{\alpha ,\beta}(d) &=&\sum_{k\ odd|\beta\ge e_k}W_3^{\alpha+e_k,
  \beta-e_k}(d) \ \ +
\\
\\
\\ && \sum_{\begin{array}{l}_{l\ge0}\\_{s\in \mathcal S_3(l,\alpha,\beta)}\end{array}} \left[
  \frac{1}{\sigma(s)} W_2\left(d-\sum d_i,0\right) 
\left(\begin{array}{c}
\frac{3d-|\alpha|+|\beta|}{2}-1
\\\frac{3d_1-|\alpha_1|+|\beta_1|}{2} ,\ldots,\frac{3d_l-|\alpha_l|+|\beta_l| }{2}
\end{array} \right)\right.
\\
\\ &&\left.\left(\begin{array}{c}
\alpha
\\  \alpha_1,\ldots, \alpha_l
\end{array} \right) 
\prod_{i=1}^l (\beta_i)_{k_i} W_3^{\alpha_i,\beta_i}(d_i) \right]

\end{array}
$$
\vspace{2ex}

where the
 integer number $\sigma(s)$ is defined as in section \ref{ch}.

Using \cite[Theorem 2]{Br7} and the technic of this paper, one
can easily prove the following result.
\begin{thm}
For any $d\ge 1$, one has $W_3(d)= (-1)^{\frac{(d-1)(d-2)}{2}}W_3^{0,de_1}(d)$.

\end{thm}

\subsection{Tropical relative Welschinger invariants}\label{trop rel inv} 
\begin{figure}[h]
\centering
\begin{tabular}{cccccc}
\includegraphics[height=4.4cm, angle=0]{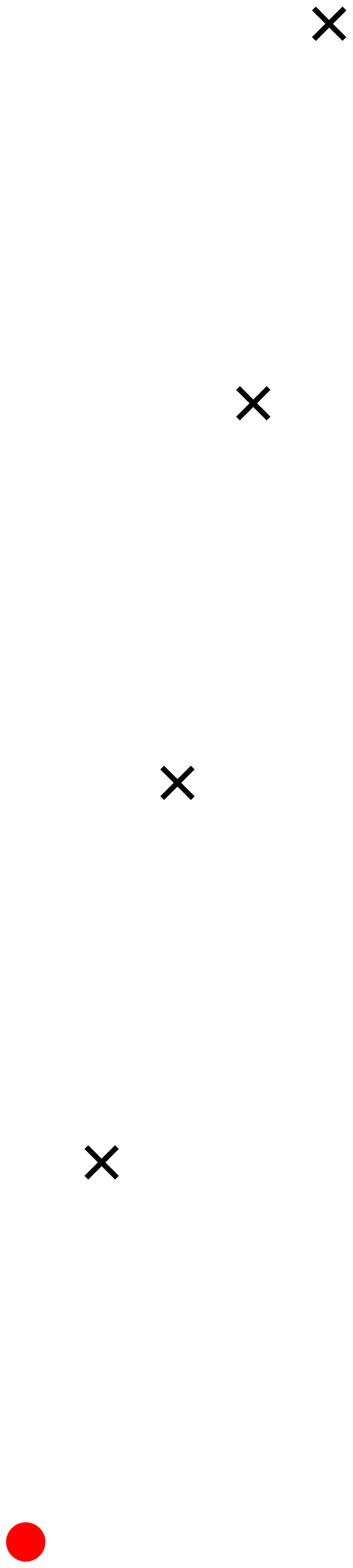}&
\includegraphics[height=4.4cm, angle=0]{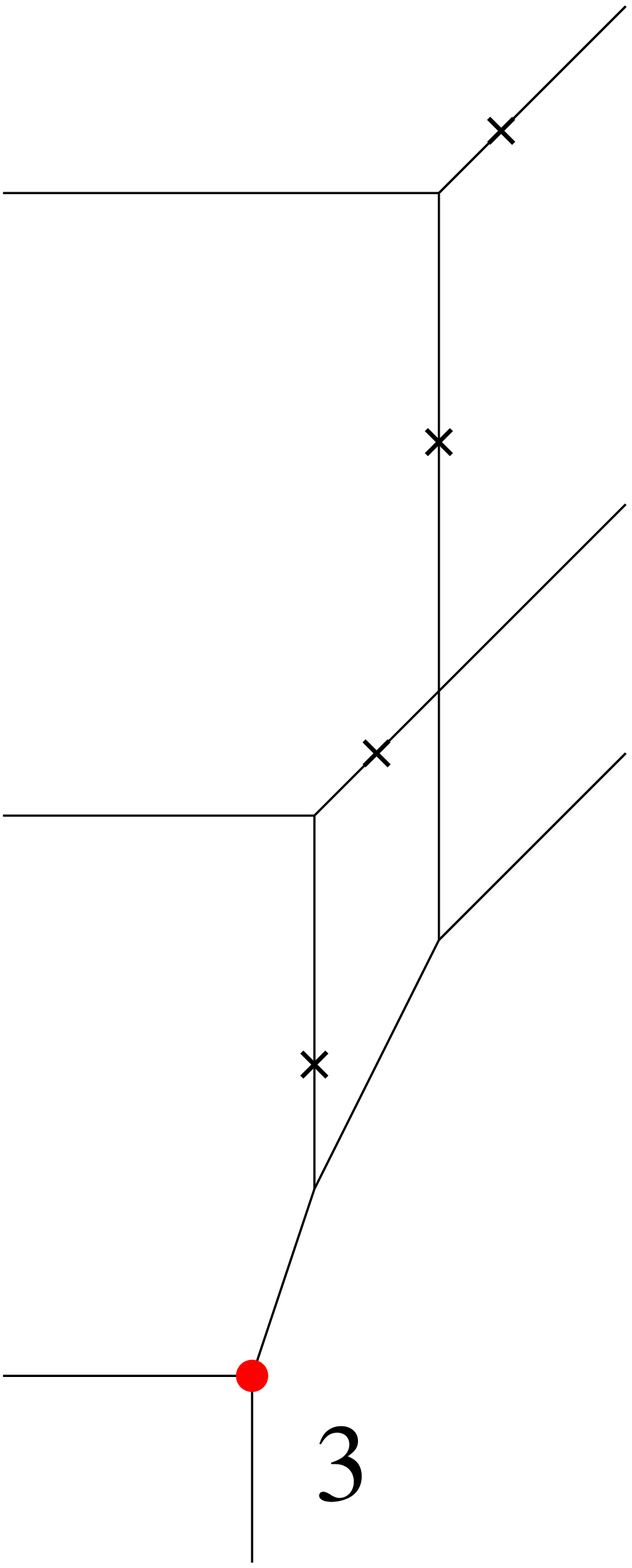}&
\includegraphics[height=4.4cm, angle=0]{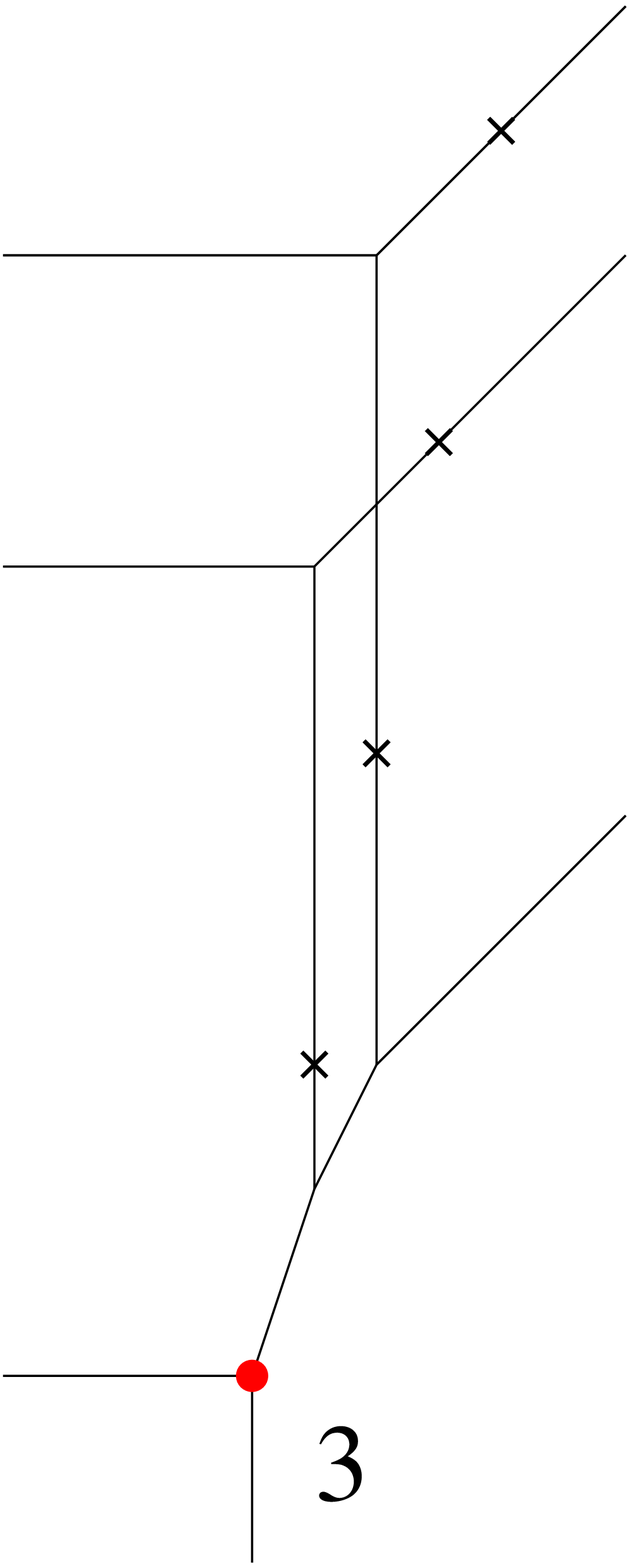}&
\includegraphics[height=4.4cm, angle=0]{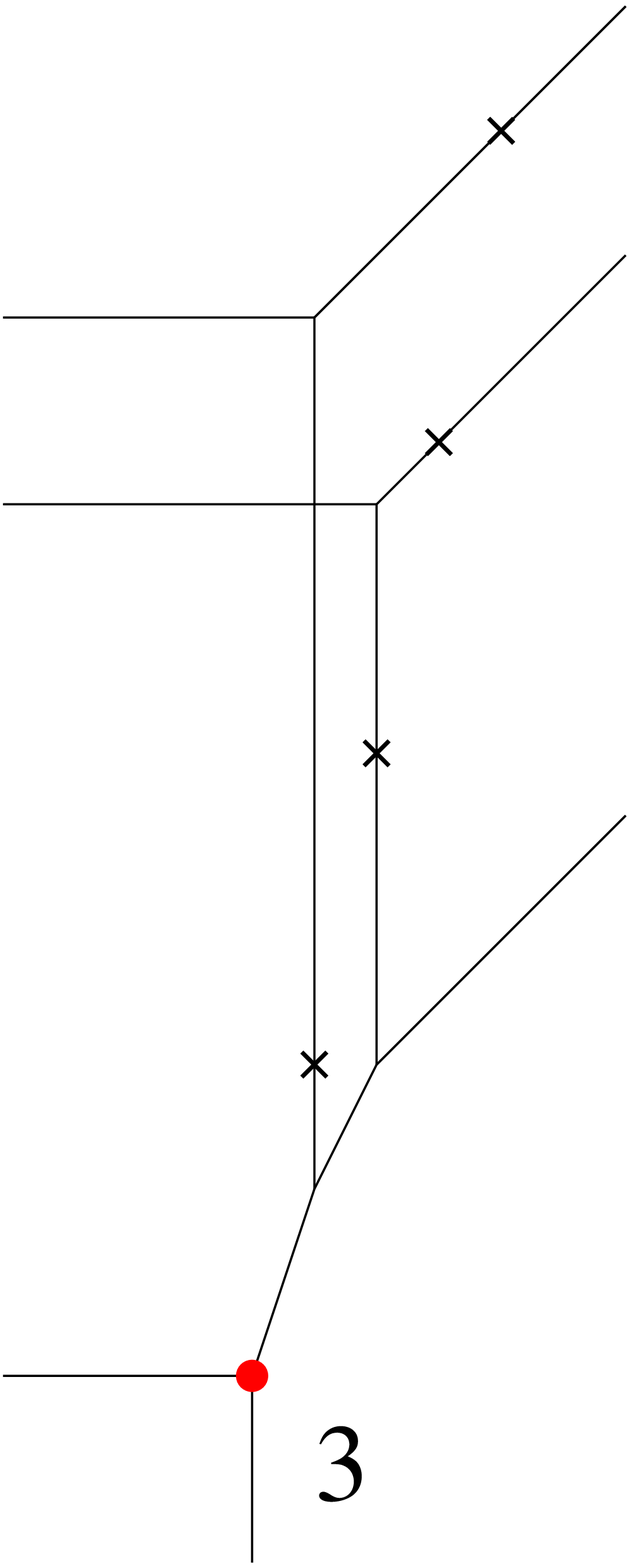}&
\includegraphics[height=4.4cm, angle=0]{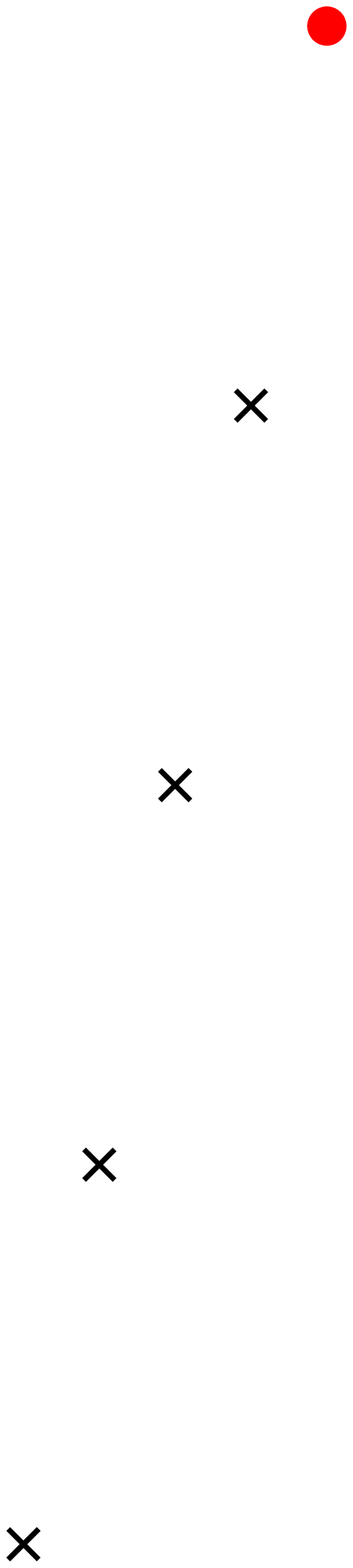}&
\includegraphics[height=4.4cm, angle=0]{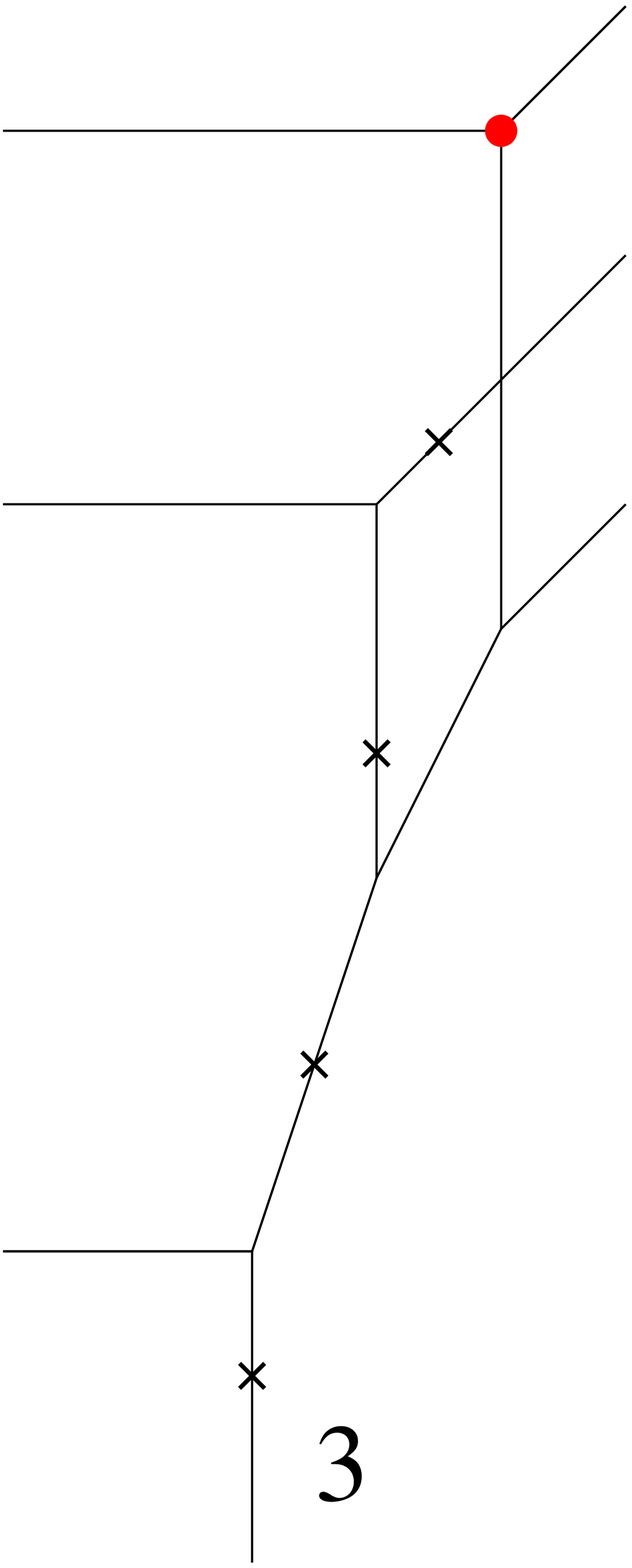}
\\
\\ a) & b) &c) &d) & e) & f)
\end{tabular}
\caption{Tropical Welschinger multiplicity does not provide tropical
  relative invariants}
\label{counter}
\end{figure}
As mentioned in section \ref{real}, one does not 
get
any relative
Welschinger invariant just by counting with Welschinger sign 
the real
rational curves of a given degree satisfying 
incidence
and tangency 
conditions.
 Surprisingly enough, Itenberg, Kharlamov and Shustin proved in
\cite{IKS3} that the situation is totally different in tropical geometry:
counting tropical curves with the tropical Welschinger multiplicity
provides tropical relative invariants in any genus.

Such a result does not hold when one considers situations related to
configurations 
of points with pairs of complex conjugated points. A counter-example is given
in Figure \ref{counter} where we look at
rational tropical cubics in $\RR^2$ with an unbounded edge of
weight 3 and passing through 5 points $x'_1,x'_2,x'_3, x'_4$ and
$x''_1$ (for brevity, we use the same notation as in \cite[Theorem
  3.1]{Sh8} and we refer to this article for more details). In Figures
\ref{counter}a and e, crosses represent the 
points $x'_i$
and the disk represents the point   $x''_1$.
There exist exactly 3 tropical cubics of genus 0 with an unbounded edge
of weight 3 passing through the configuration of  Figure 
\ref{counter}a with a non zero real multiplicity (see Figure
\ref{counter}b, c and d), and all of them are
of multiplicity 3. 
There exists only 1 tropical cubic of genus 0 with an unbounded edge
of weight 3 passing through the configuration of  Figure 
\ref{counter}e with a non zero real multiplicity (see Figure
\ref{counter}f), and it is
of multiplicity 1.

\subsection{Some computations}

We list below the first values of $N^{0,de_1}(d)$, $W_2(d,r)$ and
$W_3(d)$. 
These computations have been made using a Maple program, available at 

\noindent http://www.matcuer.unam.mx/$\sim$aubin/wel/.

$$\begin{array}{|c|c|c|c|c|c|c|c|c|c|c|}
\hline d & 1&2&3&4&5 &6&7&8&9
\\ \hline
N^{0,de_1}(d) & 1 &1  & 12&
620
&  87304 &26312976 &14616808192& 13525751027392&19385778269260800
\\\hline
\end{array}$$


\vspace{4ex}

$$\begin{array}{|c|c|c|c|c|c|}
\hline r & 0&1&2&3&4
\\ \hline
W_2(3,r)&  8& 6&4&2&0
\\\hline
\end{array}$$

\vspace{4ex}

$$\begin{array}{|c|c|c|c|c|c|c|}
\hline r & 0&1&2&3&4&5
\\ \hline
W_2(4,r)&  240& 144&80&40&16&0
\\\hline
\end{array}$$

\vspace{4ex}

$$\begin{array}{|c|c|c|c|c|c|c|c|c|}
\hline r & 0&1&2&3&4&5&6&7
\\ \hline
W_2(5,r)&  18264& 9096& 4272 &1872  &  744&248  &64  &  64
\\\hline
\end{array}$$

\vspace{4ex}

$$\begin{array}{|c|c|c|c|c|c|c|c|c|c|}
\hline r & 0&1&2&3&4&5&6&7&8
\\ \hline
W_2(6,r)&  2845440 & 1209600& 490368 &188544  &  67968&22400  &6400  &
1536 & 1024
\\\hline
\end{array}$$

\vspace{4ex}

$$\begin{array}{|c|c|c|c|c|c|c|c|c|}
\hline r & 0&1&2&3&4&5&6&7
\\ \hline
W_2(7,r)& 792731520  &293758272 & 104600448 &35670576  &11579712  &3538080  &995904  &
248976  
\\\hline
\end{array}$$

$$\begin{array}{|c|c|c|c|}
\hline r & 8&9&10
\\ \hline
W_2(7,r)& 54272 & 11776  &   -14336
\\\hline
\end{array}$$

\vspace{4ex}

$$\begin{array}{|c|c|c|c|c|c|c|c|}
\hline r & 0&1&2&3&4&5
\\ \hline
W_2(8,r)& 359935488000  &118173265920  &  37486448640  &  11463469056  &
  3367084032  & 944056320  
\\\hline
\end{array}$$

$$\begin{array}{|c|c|c|c|c|c|c|}
\hline r &6& 7& 8&9&10 & 11
\\ \hline
W_2(8,r)&249999360 & 61424640 & 13643776  &2705408    &499712 &-280576
\\\hline
\end{array}$$

\vspace{4ex}

$$\begin{array}{|c|c|c|c|c|c|c|}
\hline r & 0&1&2&3&4
\\ \hline
W_2(9,r)& 248962406889600 &73359212457600   &  20972001869568&
5807486276352
&   1553952238848
\\\hline
\end{array}$$

$$\begin{array}{|c|c|c|c|c|c|c|}
\hline r &5& 6& 7& 8&9&10 
\\ \hline
W_2(9,r)& 400246421760&   98632018560 &   23031485568&  
5021757312&  1003137408&    181785600
\\\hline
\end{array}$$

$$\begin{array}{|c|c|c|c|}
\hline r &11& 12& 13 
\\ \hline
W_2(9,r)&30391296 &   3932160 & 17326080
\\\hline
\end{array}$$

\vspace{4ex}

In particular, the values $W_2(9,12)$ and $W_2(9,13)$ disprove the
monotonicity conjecture of the 
function $r\mapsto W_2(d,r)$  by Itenberg, 
Kharlamov and Shustin (see \cite[Conjecture 6]{IKS2}). Note that the
positivity conjecture of $W_2(d,r)$ in \cite[Conjecture 6]{IKS2} has already been
disproved by Welschinger in \cite{Wel4}.

\vspace{4ex}

$$\begin{array}{|c|c|c|c|c|c|c|c|c|}
\hline d & 1&3&5&7&9 &11&13 
\\ \hline
W_3(d)& 1 &-1   & 45 & -14589 &  17756793 &-58445425017
& 426876362998821
\\\hline
\end{array}$$

\vspace{4ex}
As $W_3(2k)=0$ for any $k\ge 1$, we listed only the first values of
$W_3(2k+1)$.  

\subsection{Congruences}

Mikhalkin observed that the numbers $W_2(d,0)$ and $N_2(d)$ are equal
modulo 4. More generally, it seems that the numbers $W_2(d,r)$ and
$N_2(d)$ are equal
modulo  $2^{f(d)}$, where $f:\NN\to \NN$ is some function which goes to infinity
as $d$ goes to infinity. 
Using Kontsevich formula (see \cite{KonMan1} or \cite{DiFrIt}), one
can prove immediatly the following proposition, where $[x]$ denotes the
integer part of the real number $x$.

\begin{prop}
For any $d\ge 1$, the number $N_2(d)$ is divisible by
$2^{\left[\frac{d-1}{2} \right]}$.
\end{prop}

Up to our knowledge, no such result is known in full generality for Welschinger
invariants 
$W_2(d,r)$, although computations suggest that an analog result holds
in this case.
If $W_2(d,r)$ and
$N_2(d)$ are truly equal
modulo big powers of 2, it
 would be very interesting to find a geometrical reason. 
Note that Welschinger used symplectic field theory in \cite{Wel4} to
prove that a big power of 2 divides $W_2(d,r)$ when $r$ is close to $\frac{3d-1}{2}$.

It is proved in \cite{Br6} (see also \cite{Br7}) that the number $W_3(d)$
is also equal modulo 4 to the corresponding number of complex
curves. However, unlike in the case of enumerative invariants of $\CC P^2$, a
stronger 
congruence does not seem to hold.

\small
\def\rightmark{\em Bibliography}
\addcontentsline{toc}{section}{References}

\bibliographystyle{alpha}
\bibliography{../../Biblio.bib}

\end{document}